\documentclass[11pt,leqno]{article}
\usepackage{amsmath}
\usepackage{amssymb}
\usepackage{amsfonts}
\usepackage{enumerate}

\allowdisplaybreaks

\author{Nadira {\sc Bouchemella} and Paul {\sc Raynaud de Fitte}}

\title{Weak solutions of backward stochastic differential equations 
 with continuous generator}
 \date{August 17, 2013}

\newtheorem{theo}{Theorem}[section]
\newtheorem{prop}[theo]{Proposition}
\newtheorem{cor}[theo]{Corollary}
\newtheorem{lem}[theo]{Lemma}
\newtheorem{definition}[theo]{Definition}
\newenvironment{defe}{\begin{definition}\em}{\end{definition}}
\newtheorem{remark}[theo]{Remark}
\newenvironment{rem}{\begin{remark}\em}{\end{remark}}

\newcommand\proof{\medskip\noindent {\bf Proof}\ }
\newcommand\proofof[1]{\medskip\noindent{\bf Proof of {#1}}\ }
\newcommand\Square{\fbox{\rule{0em}{.3em}\rule{.3em}{0em}} \qquad}
\newcommand\finpr{\nobreak\hfill$\Square\qquad$\allowbreak\medskip\par}

\newcommand{\wrt}{with respect to\ }
\renewcommand\iff{if and only if\ }
\newcommand\extendedsol{extended solution}

\newcommand\cF{\mathcal{F}}
\newcommand\cG{\mathcal{G}}
\newcommand\R{\mathbb{R}}

\newcommand\un[1]{\,\rlap{{\rm 1}}\kern.22em \mbox{\rm l}_{#1}} 

\newcommand\tq{;\,} 
\newcommand\CCO[1]{\left( #1 \right)}
\newcommand\accol[1]{\left\{{#1}\right\}}
\newcommand\norm[1]{\left\Vert #1 \right\Vert}
\newcommand\abs[1]{\left\vert #1 \right\vert}
\newcommand\scal[2]{\left\langle #1,#2 \right\rangle}

\newcommand\expect{\mathop{\text{\rm E}}\nolimits}
\newcommand\expectmu{\mathop{\text{\rm E}}\nolimits} 
\newcommand\espcond[2]{\expect^{{#2}}\CCO{#1}}
\newcommand\espcondF[2]{\expect^{\cF_{{#2}}}\CCO{#1}}
\newcommand\espcondFbar[2]{\expect^{\tribuu_{{#2}}}\CCO{#1}}

\newcommand\prob{\mathop{\text{\rm P}}\nolimits}

\newcommand\law[1]{ {\mathcal L}\CCO{#1} } 
\newcommand\laws[1]{{\mathcal M}^{1,+}\CCO{{#1}}}

\newcommand\jointsol{\mathcal{M}}

\newcommand\ellp[1]{\mathop{\text{\rm L}}\nolimits^{#1}}

\newcommand\linear{\mathop{\text{\rm L}}}

\newcommand\ent[1]{\lceil #1 \rceil}

\newcommand\hyp[1]{{$(H_{#1})$}} 
\newcommand\hypp[1]{{$H_{#1}$}} 
\newcommand\Hbb{\hypp{2}}
\newcommand\Ha{\hyp{1}} 
\newcommand\Hb{\hyp{2}}
\newcommand{\jzero}{{\rm(W0)}}

\newcommand\esp{{\R^d}}
\newcommand\espY{\esp}

\newcommand\espZ{\mathbb{L}}
\newcommand\esplin{\mathbb{K}}
\newcommand\espX{\mathbb{M}}
\newcommand\espS{\mathbb{D}} 
\newcommand\espw{{\R^m}}
\newcommand\espE{\mathbb{E}} 
\newcommand\esppolish{\mathbb{S}} 

\newcommand\consf{C_f}  

\newcommand\consp{C}

\newcommand\majorr[1]{\mathfrak{M}_{#1}} 

\newcommand\CA{{\mathfrak{a}}}
\newcommand\CB{{\mathfrak{b}}}

\newcommand\TA{\mathfrak{T}}

\newcommand\pq{q}  

\newcommand\esprob{\Omega}              
\newcommand\tribu{\cF}                  
\newcommand\esprobgen{\Sigma}
\newcommand\tribugen{\cG}
\newcommand\probgen{\mathop{\text{\rm Q}}\nolimits}

\newcommand\bor[1]{\mathcal{B}\CCO{#1}}  

\newcommand\wien{W} 
\newcommand\twien{\mbox{w}} 
\newcommand\cvarqd[2]{ \left[ {#1},{#2} \right] } 
\newcommand\cvar{\mathfrak{N}} 

\newcommand\ailleurs[1]{\underline{#1}}
\newcommand\Yu{\ailleurs{Y}}
\newcommand\Zu{\ailleurs{Z}}
\newcommand\Lu{\ailleurs{L}}
\newcommand\wienu{\ailleurs{\wien}}
\newcommand\Xu{\ailleurs{X}}
\newcommand\xiu{\ailleurs{\xi}}

\newcommand\xdiff{\zeta}

\newcommand\procgene{{K}}
\newcommand\procgenne{K} 
\newcommand\vargenh{H} 
\newcommand\vargenk{K} 
\newcommand\vargena{K} 

\newcommand\VV{V} 
\newcommand\vv{v}
\newcommand\MZ{{\widehat{\VV}}} 
\newcommand\FYZ{\Sigma} 

\newcommand\coord[1]{{[#1]}}

\newcommand\mgx{N} 
\newcommand\mgxx{\widehat\mgx} 

\newcommand\ZK{\VV} 
\newcommand\zk{\vv}

\newcommand\proj[1]{\pi_{#1}} 

\newcommand\classA{{\mathcal C}} 
\newcommand\classB{{\mathcal E}} 

\newcommand\yy{y} 
\newcommand\zz{z}

\newcommand\fonctcont{\text{\rm C}} 
\newcommand\trajGen{{\Gamma}} 
\newcommand\trajY{\fonctcont_\esp[0,T]}
\newcommand\trajYD{\mathbb{D}} 
\newcommand\trajYDs{\mathbb{D}_{\mbox{\rm\tiny S}}} 
\newcommand\trajYV{\mathbb{V}} 
\newcommand\cadlag{c\`adl\`ag}

\newcommand\trajW{\fonctcont_\espw[0,T]} 

\newcommand\trajX{\fonctcont_\espX[0,T]}

\newcommand\countabl{\mathcal{N}}

\newcommand\trajZ{{\ellp{2}_\espZ[0,T]}}
\newcommand\trajZt{\ellp{2}_{\espZ}{[0,t]}}

\newcommand\Ctribu{{\mathcal D}} 
\newcommand\Htribu{{\mathcal H}} 

\newcommand\esprobL{\Gamma}
\newcommand\Ltribu{{\mathcal G}} 

\newcommand\trajZH{\mathbb{H}}
\newcommand\trajZHw{\mathbb{H}_\sigma}

\newcommand\trajE{\fonctcont_\espE[0,T]}

\newcommand\esprobb{\underline{\esprob}}
\newcommand\tribuu{\underline{\tribu}}

\newcommand\esprobc{\widetilde{\esprob}}
\newcommand\tribuc{\widetilde{\tribu}}
\newcommand\probc{\widetilde{\prob}}

\newcommand\youngs{{\mathcal Y}}
\newcommand\dirac[1]{{\delta}_{#1}} 

\newcommand\bernoul{A}

\begin{document}
\maketitle

\begin{abstract}
This paper provides a simple approach for the consideration of quadratic BSDEs with bounded terminal
 We prove the existence of a weak solution to a backward stochastic
  differential equation (BSDE) 
$$ Y_t=\xi+\int_t^T f(s,X_s,Y_s,Z_s)\,ds-\int_t^T Z_s\,d\wien_s$$ 
in a finite-dimensional space, where $f(t,x,y,z)$ is affine with
respect to $z$,  
and satisfies a sublinear growth
condition and a continuity condition. 
This solution takes the form of a triplet $(Y,Z,L)$ of processes
defined on an extended probability space and satisfying 
$$ Y_t=\xi+\int_t^T f(s,X_s,Y_s,Z_s)\,ds-\int_t^T Z_s\,d\wien_s-(L_T-L_t)$$
where $L$ is a martingale with possible jumps which is orthogonal 
to $\wien$. 
The solution is constructed on an extended probability space, 
using Young measures on the space of trajectories. 
One component of this space is 
the Skorokhod space $\trajYD$ endowed with
the topology S of Jakubowski. 
\end{abstract}


\section{Introduction}

\paragraph{Aim of the paper}
Let $(\esprob,\cF,(\cF_t)_{t\in \left[0,T \right]  },\prob)$ be a
complete probability space, 
 where $(\cF_t)_{t\geq 0}$ is the natural filtration of a standard
 Brownian motion $\wien=(\wien_t)_{t\in \left[0,T \right] }$ on 
$\espw$ and $\tribu=\tribu_T$. 

In this paper, we prove the existence of a weak solution 
(more precisely, a solution defined on an extended probability space)
to the equation
\begin{align}
 Y_t&=\xi+\int_t^T f(s,X_s,Y_s,Z_s)\,ds-\int_t^T Z_s\,d\wien_s
     -(L_T-L_t)\label{eq:BSDE-L}
\end{align} 
where $f(t,x,y,z)$ is affine with
respect to $z$,  
and satisfies a sublinear growth
condition and a continuity condition, 
$\wien$ is an $\espw$-valued standard Brownian motion,
$Y$ and $Z$ and $L$ are unknown processes, 
$Y$ and $L$ take their values in $\espY$, 
$Z$ takes its values in the space $\espZ$ of linear mappings from
$\espw$ to $\esp$, 
$\xi$ $\in\ellp{2}_\esp$ is the terminal condition, 
and 
$L$ is a 
martingale orthogonal 
to $\wien$, 
with $L_0=0$ and with \cadlag\ trajectories 
(i.e.~right continuous trajectories with left limits at every point).  
 The  process $X=(X_t)_{0\leq t\leq T}$ is $(\cF_t)$-adapted and
 continuous with values in a separable metric space  $\espX$. 
This process represents the random part of the generator $f$ and plays
a very small role in our construction. The
 space $\espX$ can be, for example, some space of trajectories, and 
$X_t$ can be, for example, the history until time $t$ of some process
$\xdiff$, 
i.e.~
$X_t=(\xdiff_{s\wedge t})_{0\leq s\leq T}$. 

Such a weak solution to \eqref{eq:BSDE-L} can be considered as
a generalized weak solution to the more classical equation
\begin{align}
 Y_t&=\xi+\int_t^T f(s,X_s,Y_s,Z_s)\,ds-\int_t^T Z_s\,d\wien_s. 
                                                    \label{eq:BSDE-gene}
\end{align}

\paragraph{Historical comments}
Existence and uniqueness of the solution $(Y,Z)$ to a 
nonlinear BSDE of the form 
\begin{align*}
 Y_t&=\xi+\int_t^T f(s,Y_s,Z_s)\,ds-\int_t^T Z_s\,d\wien_s 
\end{align*}
have been proved in
the seminal paper 
\cite{Pardoux-Peng90} by E.~Pardoux and S.~Peng, in the case when the
generator $f$ is random with
$f(.,0,0)\in\ellp{2}(\esprob\times[0,T])$, and 
$f(t,y,z)$ is Lipschitz with respect
to $(y,z)$, uniformly in the other variables. 
In \cite{Lepeltier-SanMartin}, J.P.~Lepeltier and J.~San Mart{\'\i}n  
proved in the one dimensional case the 
existence of a solution when $f$ is random, continuous with respect to $(y,z)$
and satisfies a linear growth condition 
$\norm{f(t,y,z)}\leq C(1+\norm{y}+\norm{z})$.

Equations of the form \eqref{eq:BSDE-gene},
with $f$ depending on some other process $X$,
appear in forward-backward stochastic differential equations
(FBSDEs), where $X$ is a solution of a (forward) stochastic
differential equation. 

As in the case of stochastic differential equations, one might expect
that BSDEs with continuous generator always admit at least a 
{\em weak solution},  
that is, a solution defined on a different probability space 
(generally with a larger filtration than the original one).
A work in this direction but for forward-backward stochastic
differential equations (FBSDEs) is that of K.~Bahlali,
B.~Mezerdi, M.~N'zi and Y.~Ouknine 
\cite{bahlali-mezerdi-nzi-ouknine07},
where the
original probability is changed using Girsanov's theorem. 
Let us also mention the works on weak solutions to FBSDEs 
by Antonelli and Ma
\cite{Antonelli-Ma03weak_solutions_of_FBSDEs},  
and Delarue and Guatteri 
\cite{Delarue-Guatteri06weak}, 
where the change of probability space comes from the
construction of the forward component. 

Weak solutions where the filtration is enlarged  have been studied by
R.~Buckdahn, H.J.~Engelbert and A.~R{\u{a}}{\c{s}}canu in
\cite{Buckdahn-Engelbert-Rascanu04} (see also  
\cite{Buckdahn-Engelbert05BSDE_without_strong_solution,
Buckdahn-Engelbert06continuity_of_weak}), 
using pseudopaths and the
Meyer-Zheng topology \cite{Meyer-Zheng84tightness}.  
Pseudopaths were invented by Dellacherie and Meyer 
\cite{dellacheriemeyer75book}, actually they are Young measures on the
state space (see Subsection \ref{subsect:weak_solution} for
the definition of Young measures). The success of Meyer-Zheng topology
comes from a tightness criterion which is easily satisfied and 
ensures that all limits have
their trajectories in the Skorokhod space $\trajYD$. 
We use here the fact that Meyer-Zheng's criterion also
yields tightness for Jakubowski's
stronger topology S on $\trajYD$ \cite{jakubowski97skor}. 
Note that the result of Buckdahn, Engelbert and R{\u{a}}{\c{s}}canu 
\cite[Theorem 4.6]{Buckdahn-Engelbert-Rascanu04}
is more general than ours in the sense that 
$f$ in \cite{Buckdahn-Engelbert-Rascanu04}
depends functionally 
on $Y$, more precisely, their generator $f(t,x,y)$ is defined on
 $[0,T]\times\espS\times\espS$. 
Furthermore, in \cite{Buckdahn-Engelbert-Rascanu04}, 
 $\wien$ is only supposed to be a c\`adl\`ag 
martingale. On the other hand, it is assumed in
\cite{Buckdahn-Engelbert-Rascanu04} that $f$ is bounded and does not
depend on $Z$ (but possibly on the martingale $\wien$). 
In the present paper, $f$  satisfies only a linear growth condition, 
but the main novelty (and the main difficulty) 
is that $f$ depends (linearly) on $Z$.
As our final setup is not Brownian, the process $Z$ we
construct is not directly obtained by the martingale representation
theorem, but as a limit of processes
$Z^{(n)}$ which are obtained from the martingale representation theorem.

The existence of the orthogonal component $L$ in our work comes from
the fact that our approximating sequence $(Z^{(n)})$ does not
converge in $\ellp{2}$: 
Actually it converges to $Z$ only in distribution in 
$\trajZ$ endowed with its weak topology, 
thus the stochastic integrals
$\int_0^tZ^{(n)}\,d\wien_s$ need not converge in distribution to 
$\int_0^tZ\,d\wien_s$. 
Let us mention here the work of 
Ma, Zhang and Zheng \cite{Ma-Zhang-Zheng08weak_solutions_for_FBSDEs},
on the much more intricate problem of 
existence and uniqueness of weak solutions (in the classical sense) for
forward-backward stochastic differential equations. 
Among other results, 
they prove existence of weak solutions 
with different methods and hypothesis 
(in particular the generator is assumed to be uniformly
continuous in the space of variables)  
which ensure that the approximating
sequence $Z^{(n)}$ constructed in their paper converges in $\ellp{2}$ to $Z$. 

Let us also mention the recent paper \cite{bahlali-gherbal-mezerdi11}
on the existence of an optimal control for a FBSDE. 
This optimal control and the corresponding solutions are 
obtained by taking weak limits of minimizing
controls and the corresponding strong solutions. 
The limit BSDE with the optimal control also contains an orthogonal
martingale component similar to ours.

In the case where the Brownian filtration needs to be enlarged, 
weak solutions 
are solutions which cannot be constructed as functionals of the sole Brownian
motion $\wien$. 
It is natural for this construction to add some randomness to
$\wien$ by considering Young measures on the space of trajectories of
the 
solutions we want to construct  
(let us denote momentarily $\trajGen$ this space), 
i.e.~random measures $\omega \mapsto\mu_\omega$ on $\trajGen$ 
which depend in a measurable way on the Brownian motion. The weak
solution is then constructed in the extended probability space 
$\esprobb=\esprob\times\trajGen$ with the probability 
$\mu_\omega \otimes d\prob(\omega)$. 
Young measures have been invented many times under different names for
different purposes. 
In the case of the construction of weak solutions of SDEs with
trajectories in the Skorokhod space $\trajYD$, 
they have been 
(re-)invented by Pellaumail \cite{pellaumail81weak} under the name of
{rules}. 
In the present paper, we also
construct a weak solution with the help of Young measures on a suitable
space of trajectories. 

\paragraph{Organization of the paper}
In Section \ref{sect:definitions},
we give the main definitions and hypothesis, 
in particular we discuss and compare possible definitions of weak
solutions. Using the techniques of T.G.~Kurtz, we also give 
a Yamada-Watanabe-Engelbert
type result on 
pathwise uniqueness and
existence of strong solutions. 

Section \ref{sect:construction}  
is devoted to the main result, that is, the construction of a weak solution: 
First, we construct a
sequence $(Y^{(n)},Z^{(n)})$ of strong solutions to approximating
BSDEs using a Tonelli type scheme (Subsection \ref{subsect:approximating}),
then we prove  
uniform boudedness in $\ellp{2}$ of these solutions 
(Subsection \ref{subsect:boundedness}) 
and
compactness properties in the spaces of trajectories 
(Subsection \ref{subsect:compactness}). Here the space of
trajectories is $\trajYD_\espY([0,T])\times\trajZ\times\trajYD_\espY([0,T])$, 
where $\trajYD_\espY([0,T])$ is
endowed with Jakubowski's topology S and $\trajZ$ with its weak
topology. 
Finally, 
we obtain the solution by passing to the limit of an extracted
sequence, in Young measures topology (Subsection \ref{subsect:weak_solution}). 
The proof of the main result, Theorem \ref{theo:main}, is completed 
in Subsection \ref{subsect:proof-main}.

\section{General setting, weak and strong solutions}
\label{sect:definitions}

\subsection{Generalities, 
equivalent definitions of weak solutions}
\paragraph{Notations and hypothesis}
For any separable metric space $\espE$, we denote by
$\trajE$ (respectively $\trajYD_\espE[0,T]$) 
the space of continuous (resp.~\cadlag) mappings
 on $[0,T]$ with values in $\espE$. 
(The space $\trajYD_\espY[0,T]$ will sometimes be denoted for short by
$\trajYD$.) 
Similarly, for any $q\geq 1$, 
if $\espE$ is a Banach space, and if $(\esprobgen,\tribugen,\probgen)$
is a measure space,  we denote by 
${\ellp{q}_\espE(\esprobgen)}$ the Banach space of measurable mappings 
$\varphi :\,\esprobgen\rightarrow\espE$  
such that  
$\norm{\varphi}^q_{\ellp{q}_\espE} 
:= \int_0^T \norm{\varphi(s)}^q \,d\probgen(s) < +\infty$. 

The law of a random element $X$ of a topological space  $\espE$ is
denoted by $\law{X}$. 
The conditional expectation of $X$ with respect to a $\sigma$-algebra
$\mathcal{G}$,
if it exists, 
is denoted by $\espcond{X}{\mathcal{G}}$. 
The indicator function of a set $A$ is denoted by $\un{A}$.

In the sequel, 
we are given a stochastic basis
$(\esprob,\cF,(\cF_t)_{t\in \left[0,T \right]  },\prob)$, 
the filtration $(\cF_t)$ is the filtration generated by an
$\espw$-valued standard Brownian motion $\wien$, 
augmented with the $\prob$-negligible sets, 
and $\tribu=\tribu_T$. 
We are also given 
an $\espY$-valued random variable 
$\xi$ $\in\ellp{2}_\espY(\esprob,\cF,\prob)$ (the terminal condition). 
The space of linear mappings from
$\espw$ to $\espY$ is denoted by $\espZ$. 
We denote by $\espX$ a separable metric space and by $X$ a given
$(\cF_t)$-adapted $\espX$-valued continuous process. 
Finally we are given 
a measurable mapping 
$f :\,[0,T]\times\espX\times\espY\times\espZ\rightarrow \espY$ which
satisfies 
 the following growth and continuity 
 conditions \Ha\ and \Hb\ 
(which will be needed only in Section \ref{sect:construction} for the
construction of a solution):
\begin{itemize}
  
\item[\Ha] There exists a constant $\consf\geq 0$ such that 
$\forall (t,x,y,z)\in[0,T]\times\espX\times\espY\times\espZ$, 
$\norm{f(t,x,y,z)}\leq \consf(1+\norm{z})$.

\item[\Hb]
       $f(t,x,y,z)$ is continuous \wrt $(x,y)$ and affine \wrt $z$. 
\end{itemize}

\paragraph{Weak and strong solutions}
\begin{defe}\label{def:strong}%
 A {\em strong solution} to \eqref{eq:BSDE-gene} is an $(\tribu_t)$-adapted,
$\espY\times\espZ$-valued process $(Y,Z)$ (defined on
$\esprob\times[0,T]$) 
which satisfies
\begin{align}
&\int_0^T \norm{Z_s}^2\,ds<\infty\,\prob\text{-a.e.}\label{eq:Z-L2}\\
&\int_0^T \norm{f(s,X_s,Y_s,Z_s)}\,ds<\infty\,\prob\text{-a.e.}\label{eq:f-L1}
\end{align}
and such that the BSDE \eqref{eq:BSDE-gene} holds true. 
\end{defe}
\begin{rem}\label{rem:strongL}
Similarly, a strong solution to \eqref{eq:BSDE-L} should be a triplet
$(Y,Z,L)$ defined on $\esprob\times[0,T]$) 
satisfying \eqref{eq:Z-L2}, \eqref{eq:f-L1}, and 
\eqref{eq:BSDE-L}, and such that $L$ is a \cadlag\ 
martingale orthogonal to $\wien$
and $L_0=0$, 
but this notion coincides with 
that of a strong solution to \eqref{eq:BSDE-gene}, 
because then $L$ would be an $(\cF_t)$-martingale, hence $L=0$. 
\end{rem}

\begin{rem}\label{rem:loi_de_W-X-xi}
The process $X$ is given and $(\tribu_t)$-adapted, and the final
condition $\xi$ is given and $\tribu_T$-measurable. 
By a well known result due to Doob 
(see \cite[page 603]{Doob1953book}
or \cite[page 18]{dellacheriemeyer75book}),
there exists thus a Borel-measurable mapping 
$F :\,\trajW\rightarrow\trajX\times\espY$ such that 
$(X,\xi)=F(\wien)$ a.e.
In other words, the law $\law{\wien,X,\xi}$ of 
$(\wien,X,\xi)$ is supported by the graph of $F$. 
The fact that $X$ is $(\tribu_t)$-adapted is a property of $F$: 
it means that, for every
$t\in[0,T]$, the restriction of  
$X$ to $[0,t]$ only
depends on the restriction of $\wien$ to $[0,t]$.  
\end{rem}

We now give three equivalent definitions of a weak solution:
\begin{defe}\label{def:weak}
\rule{1em}{0em}

1) A {\em weak solution} to \eqref{eq:BSDE-L} is 
a stochastic basis $(\esprobb,\tribuu,(\tribuu_t)_{0\leq t\leq T},\mu)$ 
along with a
list $(\Yu,\Zu,\Lu,\wienu, \Xu)$ of processes defined on
$\esprobb\times [0,T]$, 
and adapted to $(\tribuu_t)$, 
and a random variable $\xiu$ defined on $\esprobb$, 
such that:
 \begin{enumerate}[(W1)]

\item 
The processes $\wienu$ and $\Xu$ are continuous with values in $\espw$
and $\espX$ respectively, 
$\xiu$ takes its values in $\espY$,
and the law of $(\wienu,\Xu,\xiu)$ on $\trajW\times\trajX\times\espY$
is that of  $(\wien,X,\xi)$.

\item \label{wienu_brownian}
$\wienu$ is a standard Brownian motion with respect to the
  filtration $(\tribuu_t)$.

\item \label{traject+Lu_orthog}
The processes 
$\Yu$ and $\Lu$ are $\espY$-valued and \cadlag,
and $\Zu$ is $\espZ$-valued,
with $\expect\int_0^T\norm{\Zu_s}^2\,ds \allowbreak< \infty$, 
the process $\Lu$ is a square integrable martingale with
$\Lu_0=0$, and $\Lu$ is orthogonal
to $\wienu$.

\item  \label{YuZuLu_weaksol}
Condition \eqref{eq:f-L1} and 
the BSDE \eqref{eq:BSDE-L} hold true, replacing $Y,Z,L,\wien,X,\xi$ by
$\Yu,\Zu,\Lu,\wienu,\Xu,\xiu$.
\end{enumerate}
We then say that $(\Yu,\Zu,\Lu,\wienu, \Xu,\xiu)$ is a 
{\em weak solution defined
on $(\esprobb,\tribuu,(\tribuu_t),\mu)$}.

\smallskip
2) Following the terminology of 
\cite{jacod80weak_and_strong_solutions,Engelbert91Yamada-Watanabe, 
Kurtz07Yamada-Watanabe-Engelbert}, and   
with the preceding notations, 
the probability measure $\law{\Yu,\Zu,\Lu,\wienu, \Xu,\xiu}$ 
on $\trajYD_\espY{[0,T]}\times\trajZ\times\trajYD_\espY{[0,T]}
\times\trajW\times\trajX\times\espY$
is called a 
{\em joint solution measure 
to \eqref{eq:BSDE-L} 
generated by $(\esprobb,\tribuu,(\tribuu_t)_t,\mu)$ 
and $(\Yu,\Zu,\Lu,\wienu, \Xu,\xiu)$}. 
(Here the Borel subsets of $\trajYD_\espY{[0,T]}$ 
are generated by the projection
mappings $\proj{t} :\, x\mapsto x(t)$ for $t\in[0,T]$; we shall see
later that these sets are the Borel sets of the topology S of 
A.~Jakubowski \cite{jakubowski97skor}.) 

\smallskip
3) An {\em \extendedsol} to \eqref{eq:BSDE-L} 
consists of  
a stochastic basis $(\esprobb,\tribuu,(\tribuu_t)_{0\leq t\leq T},\mu)$ 
along with a
triplet $(Y,Z,L)$ of processes defined on $\esprobb$ such that:
 \begin{enumerate}[(E1)]
\item There exists a measurable space $(\esprobL,\Ltribu)$, and a
  filtration $(\Ltribu_t)$ on $(\esprobL,\Ltribu)$  
such that 
$$\esprobb=\esprob\times\esprobL,\quad 
\tribuu=\tribu\otimes\Ltribu,\quad
\tribuu_t=\tribu_t\otimes\Ltribu_t \text{ for every }t,
$$
and there exists a probability measure $\mu$ on $(\esprobb,\tribuu)$
such that $\mu(A\times\esprobL)=\prob(A)$ for every $A\in\tribu$.

Note that every random variable $\zeta$ defined on $\esprob$ can then be
identified with a random variable defined on $\esprobb$, by setting 
$\zeta(\omega,\gamma)=\zeta(\omega)$. Furthermore, $\tribu$ can be
viewed as a sub-$\sigma$-algebra of $\tribuu$ by identifying each
$A\in\tribu$ with the set $A\times\esprobL$. Similarly, each $\tribu_t$
can be considered as a  sub-$\sigma$-algebra of $\tribuu_t$. 
We say that 
$(\esprobb,\tribuu,(\tribuu_t)_t,\mu)$ is an {\em extension} of
  $(\esprob,\tribu,(\tribu_t)_t,\prob)$.

\item 
The process $(\wien_t)_{0\leq t\leq T}$ is 
a Brownian motion on 
 $(\esprobb,\tribuu,(\tribuu_t)_t,\mu)$
(where $\wien(\omega,\gamma):=\wien(\omega)$ for all 
$(\omega,\gamma)\in\esprobb$),

\item
The processes $Y$, $Z$ and $L$  are 
$(\tribuu_t)$-adapted,  
$Y$ and $L$ are  $\espY$-valued and \cadlag,
and $Z$ is $\espZ$-valued,
with $\expect\int_0^T\norm{Z_s}^2\,ds \allowbreak< \infty$, 
and  $L$ is a square integrable martingale with
$L_0=0$, and $L$ is orthogonal 
to $\wien$.

\item
Condition \eqref{eq:f-L1} and 
the BSDE \eqref{eq:BSDE-L} hold true.
\end{enumerate}

\end{defe}

Obviously, an \extendedsol\ is a weak solution, and a weak solution
generates a joint solution measure. Actually, these 
concepts are equivalent in the sense that: 
\begin{prop}\label{prop:equivalent_defs}
Given a joint solution measure $\nu$ to \eqref{eq:BSDE-L}, 
there exists 
an \extendedsol\ 
to \eqref{eq:BSDE-L}
which generates $\nu$. 
\end{prop}
Before we give the proof of Proposition \ref{prop:equivalent_defs},
let us give an intrinsic characterization of joint solution measures. 
Let us first observe that:

\medskip
1. It is easy to check
(see the proof of Lemma \ref{lem:equivalent-formulations})
that, if  $Y,Z,L,\wien, X,\xi$ are defined on a 
stochastic basis $(\esprobb,\tribuu,(\tribuu_t),\mu)$, 
then  \eqref{eq:BSDE-L} is equivalent to 
\begin{align}
  Y_t&=\espcondFbar{\xi+\int_t^T
    f(s,X_s,Y_s,Z_s)\,ds}{t} \label{eq:BSDE-Ftbar}\\
  \int_0^t Z_s\,d\wien_s+L_t
     &=\espcondFbar{\xi+\int_0^T f(s,X_s,Y_s,Z_s)\,ds}{t}\notag\\
     &\phantom{=f(s,X_s,Y_s^{(n)},\widetilde{Z}_s^{(n)})}
       -\espcondFbar{\xi+\int_0^T f(s,X_s,Y_s,Z_s)\,ds}{0}. \label{eq:BSDE-ZL}
\end{align}

\medskip
2. 
If  $(\Yu,\Zu,\Lu,\wienu, \Xu,\xiu)$ is a 
{weak solution defined
on a stochastic basis $(\esprobb,\tribuu,(\tribuu_t),\mu)$}, it is
still a weak solution if we reduce the filtration $(\tribuu_t)$ to a
filtration $(\tribuu'_t)$ such that $\tribuu'_t\subset\tribuu_t$ for
every $t$ and $(\Yu,\Zu,\Lu,\wienu, \Xu)$ remains adapted to
$(\tribuu'_t)$. Furthermore, conditions (W1) to (W4) remain unchanged
if we augment $(\tribuu'_t)$ with the $\mu$-negligible sets. 
So,  
if $\nu$ is a joint solution measure, there exists 
$(\Yu,\Zu,\Lu,\wienu, \Xu,\xiu)$ defined on a stochastic basis 
$(\esprobb,\tribuu,(\tribuu_t)_{0\leq t\leq T},\mu)$ 
such that 
\begin{itemize}
\item[\jzero] 
$\nu=\law{\Yu,\Zu,\Lu,\wienu, \Xu,\xiu}$
and $(\tribuu_t)=(\tribu^{\Yu,\Zu,\Lu,\wienu}_t)$, where 
$ (\tribu^{\Yu,\Zu,\Lu,\wienu}_t)$ is the filtration generated
by $(\Yu,\Zu,\Lu,\wienu)$, augmented with the $\mu$-negligible sets.
\end{itemize}

Now, Condition (W1) is clearly a condition on $\nu$. 
Let us rewrite
Conditions (W2)-(W4), under Assumption \jzero\ on 
$(\Yu,\Zu,\Lu,\wienu, \Xu,\xiu)$ and  
$(\esprobb,\tribuu,\allowbreak(\tribuu_t)_{0\leq t\leq T},\mu)$ . 
We use here techniques of Kurtz \cite{Kurtz07Yamada-Watanabe-Engelbert}. 

\medskip
$\bullet$\ 
By L\'evy's characterization of Brownian motion, 
(W\ref{wienu_brownian}) is satisfied \iff $\wienu$ is an
$(\tribuu_t)$-martingale and 
 $\cvarqd{\wienu^{\coord{i}}}{\wienu^{\coord{j}}}_t=\delta_{ij}t$ for
 all $t\in[0,T]$, where
 $\delta_{ij}$ is the Kronecker symbol and $\wienu^{\coord{i}}$ is the
 $i$th coordinate of $\wienu$. The latter condition is 
satisfied if $\wienu$ has the same law as $\wien$, thus it is 
implied by (W1). We can thus replace (W\ref{wienu_brownian}) by  
\begin{itemize}
\item[(W\ref{wienu_brownian}')] $\wienu$ is an
  $(\tribu^{\Yu,\Zu,\Lu,\wienu}_t)$-martingale. 
\end{itemize}
But (W\ref{wienu_brownian}') is equivalent to 
\begin{equation}\label{eq:Wu_mg}
\expect\bigl(
\CCO{\wienu_T-\wienu_t } 
     h(\Yu_{.\wedge t},\Zu_{.\wedge t},\Lu_{.\wedge t},\wienu_{.\wedge t})
\bigr)=0
\end{equation}
for every $t\in[0,T]$ 
and every bounded Borel measurable function $h$ defined on 
$\trajYD_\espY{[0,T]}\times\trajZ\times\trajYD_\espY{[0,T]}\times\trajW$.

\medskip
$\bullet$\ 
By \jzero, $\Yu$ and
$\Lu$ have their trajectories in $\trajYD_\espY{[0,T]}$ and $\Zu$ in
$\trajZ$, and we have $\Lu_0 =0$. 
Thus we can replace (W\ref{traject+Lu_orthog}) by 
\begin{itemize}
\item[(W\ref{traject+Lu_orthog}')] 
$\Lu$ is a square integrable $(\tribu^{\Yu,\Zu,\Lu,\wienu}_t)$-martingale 
and $\Lu$ is orthogonal to $\wienu$.
\end{itemize}
The  second part of (W\ref{traject+Lu_orthog}')
means that 
$\Lu^{(i)}\wienu^{(j)}$ is a martingale for every
  $i\in\{1,\dots,d\}$ and every $j\in\{1,\dots,m\}$, where $\Lu^{(i)}$
  and $\wienu^{(j)}$ denote the coordinate processes. Thus 
 (W\ref{traject+Lu_orthog}') can be
  expressed as 
\begin{gather}
\expect\biggl({\bigl(\Lu_T-\Lu_t\bigr)
    h(\Yu_{.\wedge t},\Zu_{.\wedge t},\Lu_{.\wedge t},\wienu_{.\wedge
      t})}\biggr)=0
\label{Lu_martingale}\\
\expect\biggl(
\CCO{\Lu^{(i)}_T\wienu^{(j)}_T-\Lu^{(i)}_t\wienu^{(j)}_t } 
     h(\Yu_{.\wedge t},\Zu_{.\wedge t},\Lu_{.\wedge t},\wienu_{.\wedge t})
\biggr)=0\label{Lu_orthog_wienu}
\end{gather}
for every $t\in[0,T]$ 
and every bounded Borel measurable function $h$ defined on 
$\trajYD_\espY{[0,T]}\times\trajZ\times\trajYD_\espY{[0,T]}\times\trajW$.

\medskip
$\bullet$\ 
Under \jzero, 
Equations \eqref{eq:BSDE-Ftbar} and \eqref{eq:BSDE-ZL} 
amount to
\begin{align*}
   \Yu_t&=\espcond{\xiu+\int_t^T
    f(s,\Xu_s,\Yu_s,\Zu_s)\,ds}{\tribu_t^{(\Yu,\Zu,\Lu,\wienu)}}\\
 \int_0^t \Zu_s\,d\wienu_s+\Lu_t
     &=\espcond{\xiu+\int_t^T
    f(s,\Xu_s,\Yu_s,\Zu_s)\,ds}{\tribu_t^{(\Yu,\Zu,\Lu,\wienu)}}\notag\\
     &\phantom{=f(s,X_s,Y_s^{(n)}}
       -\espcond{\xiu+\int_t^T
    f(s,\Xu_s,\Yu_s,\Zu_s)\,ds}{\tribu_0^{(\Yu,\Zu,\Lu,\wienu)}}. 
\end{align*}
Thus Condition  (W\ref{YuZuLu_weaksol}) is equivalent to  
\begin{gather}
\int_0^T \norm{f(s,\Xu_s,\Yu_s,\Zu_s)}\,ds<\infty
                             \,\text{ $\prob$-a.e.}\label{eq:fu-L1}\\
  \label{eq:BSDE-FtYZL-KURTZ}
   \expect\Biggl\lgroup \Bigg(
\Yu_t-\xiu-{\int_t^T f(s,\Xu_s,\Yu_s,\Zu_s)\,ds} \Biggr)
  h(\Yu_{.\wedge t},\Zu_{.\wedge t},\Lu_{.\wedge t},\wienu_{.\wedge t})
\Biggr\rgroup=0 
\end{gather}
\begin{multline}
\expect\Biggl\lgroup \Biggl(\int_0^t \Zu_s\,d\wienu_s+\Lu_t\\
     -\xiu-\int_0^T f(s,\Xu_s,\Yu_s,\Zu_s)\,ds
       +\expect\CCO{\xiu+\int_0^T f(s,\Xu_s,\Yu_s,\Zu_s)\,ds}\Bigg)\\
     \times h(\Yu_{.\wedge t},\Zu_{.\wedge t},\Lu_{.\wedge t},\wienu_{.\wedge t})
\Biggr\rgroup=0
 \label{eq:BSDE-ZLlll}
\end{multline}
for every bounded measurable function $h$ defined on 
$\trajYD_\espY{[0,T]}\times\trajZ\times\trajYD_\espY{[0,T]}\times\trajW$.


Clearly, 
Equations
\eqref{eq:Wu_mg},
\eqref{Lu_martingale},
\eqref{Lu_orthog_wienu}, 
\eqref{eq:fu-L1},
\eqref{eq:BSDE-FtYZL-KURTZ},
and \eqref{eq:BSDE-ZLlll} 
only depend on the probability measure 
$\nu=\law{\Yu,\Zu,\Lu,\wienu, \Xu,\xiu}$. 
We have thus proved the following lemma, which is actually a
characterization of joint solution measures:

\begin{lem}\label{lem:characterization_weak_sol}
Let $\Yu,\Zu,\Lu,\wienu, \Xu$ be stochastic processes defined on a 
stochastic basis 
$(\esprobb,\tribuu,(\tribuu_t)_{0\leq t\leq T},\mu)$, 
with trajectories respectively in 
$\trajYD_\espY{[0,T]}$, $\trajZ$, $\trajYD_\espY{[0,T]}$,
$\trajW$, and $\trajX$, and let $\xiu$ be an $\espY$-valued
random variable defined on $\esprobb$. 
Assume that $(\tribuu_t)$ is the filtration generated by
$(\Yu,\Zu,\Lu,\wienu)$, possibly augmented with $\mu$-negligible sets. 
Then $(\Yu,\Zu,\Lu,\wienu, \Xu,\xiu)$ is a 
weak solution to \eqref{eq:BSDE-L} defined on
$(\esprobb,\tribuu,(\tribuu_t)_{0\leq t\leq T},\mu)$ \iff 
(W1) and 
Equations
\eqref{eq:Wu_mg},
\eqref{Lu_martingale},
\eqref{Lu_orthog_wienu}, 
\eqref{eq:fu-L1},
\eqref{eq:BSDE-FtYZL-KURTZ},
and \eqref{eq:BSDE-ZLlll}
are satisfied.
\end{lem}

\begin{cor}\label{cor:sufficient_for_weaksol}
Let $\nu$ be a joint solution measure to \eqref{eq:BSDE-L}. 
Let $(\Yu,\Zu,\Lu,\wienu, \Xu,\xiu)$ and 
$(\esprobb,\tribuu,(\tribuu_t)_{0\leq t\leq T},\mu)$ 
as in Lemma
\ref{lem:characterization_weak_sol}. 
Assume that $\law{\Yu,\Zu,\Lu,\wienu, \Xu,\xiu}=\nu$. 
Then $(\Yu,\Zu,\Lu,\wienu, \Xu,\xiu)$ is a 
weak solution to \eqref{eq:BSDE-L} defined on
$(\esprobb,\tribuu,(\tribuu_t)_{0\leq t\leq T},\mu)$. 
\end{cor}
In particular, 
if $\nu$ is a joint solution measure to \eqref{eq:BSDE-L}, 
the canonical process on the space 
$\trajYD_\espY{[0,T]}\times\trajZ\times\trajYD_\espY{[0,T]}
\times\trajW\times\trajX\times\espY$ endowed with 
the probability $\nu$ 
is a weak solution to \eqref{eq:BSDE-L}.


Before we give the proof of Proposition \ref{prop:equivalent_defs},
let us give a 
definition which will be used several times.
Let $\mu$ be a probability measure on a product 
$(\esprob\times\esprobL,\tribu\otimes\Ltribu)$
of measurable spaces 
 such that $\esprobL$ is
a Polish space (or more generally, a Radon space) and $\Ltribu$ is its
Borel $\sigma$-algebra. Let $\prob$ denote the marginal measure of $\mu$ on
$\esprob$, that is, $\prob(A)=\mu(A\times\esprobL)$ for all
$A\in\tribu$. 
Then 
there exists a unique (up to equality $\prob$-a.e.) 
family $(\mu_\omega)_{\omega\in\esprob}$ such that
$\omega\mapsto\mu_\omega(B)$ is measurable for every $B\in\Ltribu$, and 
\begin{equation}\label{eq:disinteg}
\mu(\varphi)=\int_\esprob\mu_\omega(\varphi(\omega,.))\,d\prob(\omega)
\end{equation}
for every $\tribu\otimes\Ltribu$-measurable nonnegative function 
$\varphi :\,\esprob\times\esprobL\rightarrow\R$, 
see e.g.~\cite{valadier73desi}.
\begin{defe}\label{def:disintegration} 
The family $(\mu_\omega)$ in \eqref{eq:disinteg} 
is called the {\em disintegration} of $\mu$
with respect to $\prob$. 
It is convenient to denote 
$$\mu=\mu_\omega\otimes d\prob(\omega).$$
\end{defe}

\proofof{Proposition \ref{prop:equivalent_defs}}
Let 
$\esprobL=\trajYD_\espY{[0,T]}\times\trajZ\times\trajYD_\espY{[0,T]}$.  
Let $\Ltribu$ be the Borel $\sigma$-algebra of $\esprobL$, and, for
each $t\in[0,T]$, let $\Ltribu_t$ be the $\sigma$-algebra 
generated by the projection of $\esprobL$ onto 
$\trajYD_\espY{[0,t]}\times\trajZt\times\trajYD_\espY{[0,t]}$. 
Let $F :\,\trajW\mapsto\trajX\times\espY$ 
be as in Remark \ref{rem:loi_de_W-X-xi}. 
Then, with slight abuses of notations, 
$\nu$ is the image of a probability measure $\lambda$ on 
$\esprobL\times\trajW\times\trajX\times\espY$
by the mapping
$$
\left\{\begin{array}{lcl}
\esprobL\times\trajW
&\rightarrow&
\esprobL\times\trajW \times \trajX\times\espY\\
(y,z,l,\twien)&\mapsto&(y,z,l,\twien,F(\twien)). 
\end{array}\right.
$$
Let  $(\lambda_{\tiny\twien})_{\tiny\twien\in\trajW}$ be the
disintegration of $\lambda$ with respect to $\law{\wien}$, that is, 
$(\lambda_{\tiny\twien})_{\tiny\twien\in\trajW}$ is  
a family of probability measures on $(\esprobL,\Ltribu)$
such that, for every bounded measurable 
$\varphi :\,\esprobL\times\trajW \rightarrow\R$, 
$$\lambda(\varphi)=\int_{\trajW}
  \CCO{
  \int_{\esprobL}
         \varphi{(y,z,l,\twien)}\,d\lambda_{\tiny\twien}(y,z,l)
       }
                          \,d\law{\wien}(\twien).
$$
Now, let 
$$\esprobb=\esprob\times\esprobL, 
\quad \tribuu=\tribu\otimes\Ltribu, \quad 
\tribuu_t=\tribu_t\otimes\Ltribu_t \quad (t\in[0,T]),$$  
and let $\mu=\lambda_{\wien(\omega)}\otimes d\prob(\omega)$,
i.e.~$\mu$ is  the probability measure on
$(\esprobb,\tribuu)$ 
such that 
$$\mu(\varphi)=\int_{\esprob}
  \CCO{
  \int_{\esprobL}
         \varphi{(\omega,y,z,l)}\,d\lambda_{\wien(\omega)}(y,z,l)
       }
                          \,d\prob(\omega)
$$
for every bounded measurable 
$\varphi :\,\esprobb \rightarrow\R$. 
We define the random variables $Y$, $Z$, $L$, $\wien$, $X$ and $\xi$
on $\esprobb$ by 
\begin{gather*}
Y(\omega,y,z,l)=y,\quad 
Z(\omega,y,z,l)=z,\quad
L(\omega,y,z,l)=l,\\
\wien(\omega,y,z,l)=\wien(\omega),\quad
X(\omega,y,z,l)=X(\omega),\quad 
\xi(\omega,y,z,l)=\xi(\omega).
\end{gather*}
Then $\tribuu_t=\tribu^{(Y,Z,L,\wien)}_t$ for every
$t\in[0,T]$, where 
$ (\tribu^{Y,Z,L,\wien}_t)$ is the filtration generated
by $(Y,Z,L,\wien)$ augmented with the $\prob$-negligible sets. 
Furthermore, we have 
$\law{Y,Z,L,\wien, X,\xi}=\nu$, thus, 
by Corollary \ref{cor:sufficient_for_weaksol}, 
 as $\nu$ is a joint solution
measure to \eqref{eq:BSDE-L},  
$(Y,Z, L,\wien,X,\xi)$ is a weak solution  to \eqref{eq:BSDE-L} defined on 
$(\esprobb,\tribuu,(\tribuu_t),\mu)$. 
\finpr
\begin{rem}
The extension which generates $\nu$ in Proposition \ref{prop:equivalent_defs} 
is not unique. 
The one we construct 
in Section \ref{sect:construction} 
is based on a different construction of the auxiliary space $\esprobL$. 
\end{rem}

We now give a criterion for an extended probability space to preserve
martingales. 
The equivalence 
\eqref{extendedM}$\Leftrightarrow$\eqref{verygoud} 
in the following lemma 
is contained in Lemma 2.17 of \cite{jacod-memin81weakstrong}. 
\begin{lem}\label{lem:extended_martingale}
Let $(\esprobb,\tribuu,(\tribuu_t),\mu)
=(\esprob\times\esprobL,
\tribu\otimes\Ltribu,
(\tribu_t\otimes\Ltribu_t),\mu)$ 
be an extension of $(\esprob,\tribu,(\tribu_t),\prob)$. 
Let $(\mu_\omega)$ be the {disintegration} of $\mu$ with respect to
$\prob$. 
The following are equivalent:
\begin{enumerate}[(i)]
\item \label{extendedW}
$\wien$ is an $(\tribuu_t)$-Brownian motion under $\mu$,
\item \label{extendedM}
Every  $(\tribu_t)$-martingale is an $(\tribuu_t)$-martingale
  under $\mu$\footnotemark,
  \footnotetext{According to Jacod and M\'emin's terminology 
  \cite[Definition 1.7]{jacod-memin81weakstrong}, 
  this means that $(\esprobb,\tribuu,(\tribuu_t)_t,\mu)$ 
  is a {\em very good extension} 
  of $(\esprob,\tribu,(\tribu_t)_t,\prob)$.
  A similar condition is 
  called {\em compatibility} in \cite{Kurtz07Yamada-Watanabe-Engelbert}.} 
\item \label{verygoud}
For every $t\in[0,T]$ and every $B\in\Ltribu_t$, the mapping
  $\omega\mapsto\mu_\omega(B)$ is $\tribu_t$-measurable. 
\end{enumerate}
\end{lem}
\proof
Assume \eqref{extendedW}. 
Let $M$ be an $(\tribu_t)$-martingale with values in
$\R^k$ for some integer $k$. Assume first that $M$ is  square
integrable. 
By the martingale representation theorem, 
there exists an
$(\tribu_t)$-adapted process $H$ 
with $\expect\int_0^T H_s^2\,ds<+\infty$ such that 
$M_t=M_0+\int_0^tH_s\,d\wien_s$. 
By \eqref{extendedW}, $M$ is an $(\tribuu_t)$-martingale. 
In the general case, denote, 
for every integer $N\geq 1$,  
$$M^N_T=
\begin{cases}
\frac{M_T}{\norm{M_T}}&\text{ if }\norm{M_T}>N\\
M_T&\text{ if }\norm{M_T}\leq N,
\end{cases}
$$
and set $M^N_t=\espcondF{M^N_T}{t}$ for $0\leq t\leq T$. Then, for any
$A\in\tribuu_t$, using Lebesgue's dominated convergence theorem, we
have
\begin{equation*}
\expect\CCO{\un{A}\CCO{M_T-M_t}}
=\lim_{N\rightarrow\infty}\expect\CCO{\un{A}\CCO{M^N_T-M^N_t}}=0
\end{equation*}
which proves that $M_t=\espcondFbar{M_T}{t}$. Thus
\eqref{extendedM} is satisfied.

Assume \eqref{extendedM}, and let $B\in\Ltribu_t$.  
For $u :\, \esprob\rightarrow\R$ 
and $v :\, \esprobL\rightarrow\R$ 
we denote by $u\otimes v$ the function defined on
$\esprob\times\esprobL$ by 
$u\otimes v(\omega,x)=u(\omega)v(x)$. 
For each bounded $\tribu$-measurable random variable $\vargena$, we have 
\begin{align*}
\expect\CCO{\vargena\,\mu_.(B)}
&=\mu\CCO{\vargena\otimes\un{B}}
=\mu\CCO{\espcondFbar{\vargena\otimes\un{B}}{t}}
=\mu\CCO{\espcondF{\vargena}{t}\otimes\un{B}}\notag\\
&=\expect\CCO{\espcondF{\vargena}{t}\,\mu_.(B)}
=\expect\CCO{\espcondF{\vargena}{t}\,\espcondF{\mu_.(B)}{t}}\notag\\
&=\expect\CCO{{\vargena}\,\espcondF{\mu_.(B)}{t}},
\end{align*}
which yields
$\mu_.(B)=\espcondF{\mu_.(B)}{t}$.
Thus $\mu_.(B)$ is $\tribu_t$-measurable,
which proves \eqref{verygoud}.

Assume \eqref{verygoud}. 
To prove \eqref{extendedW}, 
we only need to check that 
$\wien$ has independent increments under $\mu$. 
Let $t\in[0,T]$, and let $s>0$ such that $t+s\in[0,T]$. 
Let us prove that, for 
any $A\in\tribuu_t$ 
and any  Borel subset $C$ of $\espw$, we have 
\begin{equation}\label{eq:independence}
\mu\CCO{A\cap {\{\wien_{t+s}-\wien_{t}\in C\}}}
=\mu(A)\, \mu{\{\wien_{t+s} -\wien_{t} \in C\}}.
\end{equation} 
Let $B=\{\omega\in\esprob\tq
\wien_{t+s} {(\omega)}-\wien_{t} {(\omega)}\in C\}$.  
We have
\begin{align*}
\mu\CCO{A\cap (B\times\esprobL)}
&=\int_{\Omega\times\esprobL} 
\un{A}(\omega,\gamma))\un{B}(\omega)\,d\mu(\omega,\gamma) \\
&=\int_\Omega\mu_\omega(\un{A}(\omega,.))\un{B}(\omega)\,d\prob(\omega) \\
&=\int_\Omega\mu_\omega(\un{A}(\omega,.))\,d\prob(\omega)\,\prob(B) \\
&=\mu(A)\, \mu(B\times\esprobL),
\end{align*}
which proves \eqref{eq:independence}. Thus 
$\wien_{t+s} -\wien_{t} $ is independent of $\tribuu_t$. 
\finpr

\subsection{Pathwise uniqueness and strong solutions}
One easily sees that, under hypothesis \Ha\ and \Hb,
  Equation \eqref{eq:BSDE-gene} may have infinitely many {\em strong}
  solutions. 
For example, let $d=m=1$, $\xi=0$, and
$f(s,x,y,z)=\sqrt{\abs{y}}$. Then, for any $t_0\in[0,T]$, we get a
solution by setting $Z=0$ and 
$$Y_t=\begin{cases}
\frac{1}{4}(t_0-t)^2& \text{if } 0\leq t\leq t_0\\
0& \text{if } t_0\leq t\leq T. 
\end{cases}$$

Following the usual terminology, let us say that 
{\em pathwise uniqueness} holds for Equation \eqref{eq:BSDE-L} if two
weak solutions defined on the same probability space,  
and with respect to the same $(\wien,X,\xi)$, 
necessarily coincide. 
Thus, in our setting, pathwise uniqueness does not necessarily hold.

T.~G.~Kurtz \cite{Kurtz07Yamada-Watanabe-Engelbert}
has proved a very general version of the Yamada-Watanabe and Engelbert
theorems on uniqueness and existence of strong solutions to
stochastic equations, which includes SDEs, BSDEs and FBSDEs, but
without $z$ in the generator. 
His results are based on the convexity of the set of joint
solution-measures when the trajectories lie in a Polish space.

We can consider here that $\trajYD_\espY{[0,T]}$ is equipped with Skorokhod's 
topology $J_1$, which is Polish 
(actually, in Section \ref{sect:construction},  
 we will use Jakubowski's topology S on $\trajYD_\espY{[0,T]}$, 
which is not Polish, 
but this topology has the same Borel subsets as $J_1$). 
Thus the space 
$\esprobL=\trajYD_\espY{[0,T]}\times\trajZ\times\trajYD_\espY{[0,T]}$
is Polish.  
In particular, Theorem 3.15 of \cite{Kurtz07Yamada-Watanabe-Engelbert} 
applies to our
framework.

\begin{prop}\label{prop:yamada-watanabe}%
{\bf (Yamada-Watanabe-Engelbert \`a la Kurtz)}
Assume that 
pathwise uniqueness holds for Equation \eqref{eq:BSDE-L}. 
Then every weak solution to \eqref{eq:BSDE-L} 
is a strong solution. 
Conversely, 
if every solution to \eqref{eq:BSDE-gene} is strong, 
(equivalently, by Remark \ref{rem:strongL}, 
if every solution to \eqref{eq:BSDE-L} is strong), then 
pathwise uniqueness holds for Equation \eqref{eq:BSDE-L}. 
\end{prop}
\proof
In order to apply 
\cite[Theorem 3.15]{Kurtz07Yamada-Watanabe-Engelbert}, 
we only need to check that
the set of joint solution measures to
\eqref{eq:BSDE-L} is convex. 
(Theorem 3.15 in \cite{Kurtz07Yamada-Watanabe-Engelbert} supposes that 
$\mu\in S_{\Gamma,C,\nu}$ in the notations of \cite{Kurtz07Yamada-Watanabe-Engelbert},
but a joint solution measure  
is exactly an element of
$S_{\Gamma,C,\nu}$.) 
We check this convexity by an adaptation of 
\cite[Example 3.17]{Kurtz07Yamada-Watanabe-Engelbert}.

The set $\jointsol$ of laws of joint solution measures to
\eqref{eq:BSDE-L} 
is the set of probability laws 
of 
random variables  
$(\Yu,\Zu,\Lu,\wienu,\Xu,\xiu)$ 
with values in
$\trajYD_\espY{[0,T]}\times\trajZ\times\trajYD_\espY{[0,T]}\times\trajW
\times\trajX\times\espY$,  
satisfying the conditions of Lemma
\ref{lem:characterization_weak_sol}. But each of these conditions is a
convex constraint on $\jointsol$.

For example, to show that Equation \eqref{eq:BSDE-ZLlll} 
is a convex constraint on $\jointsol$, let us prove that the 
map $\law{Z,M}\mapsto\law{\int_0^. Z_s\,dM_s}$ preserves convex
combinations of probability laws. More precisely, 
let $\laws{\mathfrak{X}}$ denote the set of all probability
laws on a measurable space $\mathfrak{X}$. 
Let $\mathcal{C}$ be the subset of $\laws{\trajZ\times\trajW}$
consisting of laws of processes $(Z,M)$ such that 
$M$ is a standard $\espw$-valued Brownian motion 
and $Z$ is $\espZ$-valued and $M$-adapted. 
We show that the mapping 
$$
\left\{\begin{array}{lcl}
\mathcal{C}&\rightarrow&\laws{\trajY}\\
\law{Z,M}&\mapsto&\law{\int_0^. Z_s\,dM_s}
\end{array}\right.
$$ 
preserves convex
combinations of probability laws. 
Indeed, Let $\mu_1,\mu_2\in\mathcal{C}$, and let $p\in[0,1]$. 
Let $(Z^1,M^1)$ and $(Z^2,M^2)$ be adapted processes defined on
stochastic bases $(\esprob_1,\tribu^1,(\tribu^1_t),\prob_1)$ 
and $(\esprob_2,\tribu^2,(\tribu^2_t),\prob_2)$ with laws $\mu_1$ and $\mu_2$
respectively, such that $M^1$ (respectively $M^2$) is an
$(\tribu^1_t)$-Brownian motion 
(resp. $(\tribu^2_t)$-Brownian motion). Let $\bernoul$ be a random
variable taking the values $1$ with probability $p$ and $-1$ with
probability $1-p$, 
defined on a probability space $(\esprob_0,\tribu^0,\prob_0)$. 
We define a stochastic basis $(\esprobc,\tribuc,(\tribuc_t),\probc)$ by
\begin{gather*}
\esprobc=\esprob_0\times\esprob_1\times\esprob_2, \quad
\tribuc=\tribu^0\otimes\tribu^1\otimes\tribu^2,\quad
\tribuc_t=\tribu^0\otimes\tribu^1_t\otimes\tribu^2_t,\\
\probc=\prob_0\otimes\prob_1\otimes\prob_2.
\end{gather*}
For $(\omega_0,\omega_1,\omega_2)\in\esprobc$, set 
\begin{align*}
(Z,M)(\omega_0,\omega_1,\omega_2)=
\begin{cases}
(Z^1,M^1)(\omega_1) & \text{ if }\bernoul(\omega_0)=1\\
(Z^2,M^2)(\omega_2) & \text{ if }\bernoul(\omega_0)=-1.
\end{cases}
\end{align*}
Then 
$\law{Z,M}=p\mu_1+(1-p)\mu_2\in\mathcal{C}$, 
the process $M$ is $(\tribuc_t)$-Brownian, and 
$$\int_0^.Z_s\,dM_s
=\un{\{\bernoul=1\}}\int_0^.Z^1_s\,dM^1_s
+\un{\{\bernoul=-1\}}\int_0^.Z^2_s\,dM^2_s,$$
thus 
$$\law{\int_0^.Z_s\,dM_s}
=p\,\law{\int_0^.Z^1_s\,dM^1_s}+(1-p)\,\law{\int_0^.Z^2_s\,dM^2_s}.$$

The same technique can be applied to show that 
Equations
\eqref{eq:Wu_mg},
\eqref{Lu_martingale},
\eqref{Lu_orthog_wienu}, 
\eqref{eq:fu-L1},
\eqref{eq:BSDE-FtYZL-KURTZ},
and \eqref{eq:BSDE-ZLlll} are convex constraints on $\jointsol$. 
Thus $\jointsol$ is convex. 

\finpr


\section{Construction of a weak solution}
\label{sect:construction}
\begin{theo}\label{theo:main}
Assume that $f$ satisfies hypotheses \Ha\ and \Hb. 
Then Equation \eqref{eq:BSDE-L} 
admits a weak solution. 
\end{theo}

This section is entirely devoted to the proof of Theorem
\ref{theo:main}, 
by constructing an 
\extendedsol\ to \eqref{eq:BSDE-L} 
in the terminology of Definition \ref{def:weak}.

In Subsections \ref{subsect:approximating} to  
\ref{subsect:weak_solution}, 
we only assume that $f$ is measurable 
and satisfies the growth condition \Ha. 
Condition \Hb\ will be needed only 
in Subsection \ref{subsect:proof-main},
for the final part of the proof of Theorem \ref{theo:main}.

Note that the counterexample given by Buckdahn and Engelbert in 
\cite{Buckdahn-Engelbert05BSDE_without_strong_solution} 
does not fit in our framework, 
and we do not know any example of a BSDE of the form \eqref{eq:BSDE-gene} 
or \eqref{eq:BSDE-L}
under hypothesis \Ha\ and \Hb\ which has no strong solution.

\subsection{Construction of an approximating sequence 
of solutions} 
\label{subsect:approximating}


\paragraph{Approximating equations}
The proof of Lemma \ref{lem:equivalent-formulations} 
will show that
\eqref{eq:BSDE-gene} amounts to the following 
equations \eqref{eq:BSDE-Ft} and \eqref{eq:BSDE-Z}:
\begin{align}
  Y_t&=\espcondF{\xi+\int_t^T
    f(s,X_s,Y_s,Z_s)\,ds}{t} \label{eq:BSDE-Ft}\\
  \int_0^t Z_s\,d\wien_s
     &=\espcondF{\xi+\int_0^T f(s,X_s,Y_s,Z_s)\,ds}{t}\notag\\
     &\phantom{=f(s,X_s,Y_s^{(n)},\widetilde{Z}_s^{(n)})}
       -\expect\CCO{\xi+\int_0^T f(s,X_s,Y_s,Z_s)\,ds} \label{eq:BSDE-Z}.
\end{align} 
We can now write the approximating equations for \eqref{eq:BSDE-Ft}
and \eqref{eq:BSDE-Z}: 
\begin{align}
  Y_t^{(n)}&=\espcondF{\xi+\int_{t+1/n}^T
    f(s,X_s,Y_s^{(n)},\widetilde{Z}_s^{(n)})\,ds}{t} \label{eq:BSDE-Ftn}\\
  \int_0^t Z_s^{(n)}\,d\wien_s
     &=\espcondF{\xi+\int_0^T
       f(s,X_s,Y_s^{(n)},\widetilde{Z}_s^{(n)})\,ds}{t}\notag\\
     &\phantom{=f(s,X_s,Y_s^{(n)},\widetilde{Z}_s^{(n)})}
    -\expect\CCO{\xi+\int_{0}^T 
             f(s,X_s,Y_s^{(n)},\widetilde{Z}_s^{(n)})\,ds}. \label{eq:BSDE-Zn}
\end{align}
Here and in the sequel, 
\begin{itemize}
\item $f$ is extended by setting $f(t,x,y,z)=0$ for $t>T$; similarly,
  for any function or process $v$ defined on $[0,T]$, we set $v(t)=0$
  for $t>T$, 
\item we denote
  $\widetilde{Z}_s^{(n)}=\espcondF{{Z}_{s+1/n}^{(n)}}{s}$. 
\end{itemize}

\begin{prop}\label{prop:YnZn_inegrables}
The system \eqref{eq:BSDE-Ftn}-\eqref{eq:BSDE-Zn}
admits a unique strong solution $(Y^{(n)},Z^{(n)})$. 
Furthermore, for every $n\geq 1$, 
$Y^{(n)}_t \in\ellp{2}_{\esp}(\esprob)$ for each $t\in[0,T]$
and  $Z^{(n)}\in\ellp{2}_{\espZ}(\esprob\times[0,T])$. 
\end{prop}
\proof 
Let $T_k=T-\frac{k}{n}$, $k=0,\dots,\ent{nT}$, 
where $\ent{nT}$ 
is the
integer part of $nT$. 
Observe first that for each $k$, 
\eqref{eq:BSDE-Zn} amounts on the interval $\left]T_{k+1},T_k\right] $ to 
\begin{multline}
   \int_{T_{k+1}}^t Z_s^{(n)}\,d\wien_s
=\espcondF{\xi+\int_{T_{k+1}}^Tf(s,X_s,Y_s^{(n)},\widetilde{Z}_s^{(n)})\,ds}{t}\\
-\espcondF{\xi+\int_{T_{k+1}}^Tf(s,X_s,Y_s^{(n)},\widetilde{Z}_s^{(n)})\,ds}{T_{k+1}}.
\label{eq:BSDE-Znk}
\end{multline}

Now, the construction of $(Y^{(n)},Z^{(n)})$ is easy by backward induction:
For $T_1\leq t\leq T=T_0$, we have 
  $Y_t^{(n)}=\espcondF{\xi}{t}$ 
and $(Z_t^{(n)})_{T_1\leq t\leq T}$ is the unique predictable process
 such that 
$\expect\int_{T_1}^T \CCO{Z_t^{(n)}}^2\,ds <+\infty$ 
and
\begin{multline*}
  \int_{T_{1}}^t Z_s^{(n)}\,d\wien_s
=\espcondF{\xi+\int_{T_{1}}^Tf(s,X_s,Y_s^{(n)},0)\,ds}{t}\\
-\espcondF{\xi+\int_{T_{1}}^Tf(s,X_s,Y_s^{(n)},0)\,ds}{T_{1}}.
\end{multline*}
Suppose $(Y^{(n)},Z^{(n)})$ is defined on the time interval
$\left]T_k,T\right]$, with 
$k<\ent{nT}$, 
then $Y^{(n)}$ is defined in a unique way on $\left]T_{k+1},T_k\right]$
by \eqref{eq:BSDE-Ftn} and then 
$Z^{(n)}$ on the same interval by \eqref{eq:BSDE-Znk}. 
Furthermore, we get by induction 
from \eqref{eq:BSDE-Znk} 
that 
$Z^{(n)}\in\ellp{2}_{\espZ}(\esprob\times[0,T])$. Then, using this
latter result in \eqref{eq:BSDE-Ftn}, 
we deduce 
that $Y^{(n)}_t
\in\ellp{2}_{\esp}(\esprob)$ for each $t\in[0,T]$. 
\finpr

The following result links \eqref{eq:BSDE-Ftn}
and \eqref{eq:BSDE-Zn} to an approximate version of \eqref{eq:BSDE-gene}:
\begin{lem}\label{lem:equivalent-formulations}
Equations \eqref{eq:BSDE-Ftn} and \eqref{eq:BSDE-Zn} are equivalent to 
\begin{align}
   Y_t^{(n)}&=\xi+\int_{t}^T f(s,X_s,Y_s^{(n)},\widetilde{Z}_s^{(n)})\,ds
        -\int_t^T Z_s^{(n)}\,d\wien_s 
        -U_t^{(n)}\label{eq:BSDE-genen}\\
\intertext{with $Y^{(n)}$ adapted and}
 U_t^{(n)}
&=\espcondF{\int_t^{t+1/n}f(s,X_s,Y_s^{(n)},\widetilde{Z}_s^{(n)})\,ds}{t}.\notag
\end{align}
\end{lem}
\proof
Assume \eqref{eq:BSDE-Ftn} and \eqref{eq:BSDE-Zn}.  
Denoting 
$$M_t^{(n)}=\espcondF{\xi+\int_0^T f(s,X_s,Y_s^{(n)},\widetilde{Z}_s^{(n)})\,ds}{t}
=M_0^{(n)}+\int_0^t Z_s^{(n)}\,d\wien_s,$$ 
we get 
\begin{multline*}
  M_t^{(n)}=\espcondF{\xi+\int_{t+1/n}^T
    f(s,X_s,Y_s^{(n)},\widetilde{Z}_s^{(n)})\,ds}{t}\\
    +\int_0^t f(s,X_s,Y_s^{(n)},\widetilde{Z}_s^{(n)})\,ds+U_t^{(n)}.
\end{multline*}
By \eqref{eq:BSDE-Ftn}, this yields
\begin{align*}
  M_t^{(n)}
 &=Y_t^{(n)}+\int_0^t f(s,X_s,Y_s^{(n)},\widetilde{Z}_s^{(n)})\,ds+U_t^{(n)},
\end{align*}
that is,
\begin{align*}
  Y_t^{(n)}&=M_t^{(n)}-\int_0^t
  f(s,X_s,Y_s^{(n)},\widetilde{Z}_s^{(n)})\,ds-U_t^{(n)}\\ 
       & =M_0^{(n)}
       + \int_0^t Z_s^{(n)}\,d\wien_s
       -\int_0^t f(s,X_s,Y_s^{(n)},\widetilde{Z}_s^{(n)})\,ds-U_t^{(n)}.\\
\intertext{
In particular,} 
Y_T^{(n)}&=\xi=M_0^{(n)}+\int_0^T Z_s^{(n)}\,d\wien_s-\int_0^T
f(s,X_s,Y_s^{(n)},\widetilde{Z}_s^{(n)})\,ds\\
\intertext{thus}
  Y_t^{(n)}-Y_T^{(n)}
&=-\int_t^T Z_s^{(n)}\,d\wien_s+\int_t^T
              f(s,X_s,Y_s^{(n)},\widetilde{Z}_s^{(n)})\,ds-U_t^{(n)}, 
\end{align*}
which proves \eqref{eq:BSDE-genen}.

Conversely, assume \eqref{eq:BSDE-genen} and that $Y^{(n)}$ is
adapted. 
Denote ${\VV}_t^{(n)}=\int_0^t Z_s^{(n)}\,d\wien_s$. 
We have 
\begin{align*}
 Y_t^{(n)}&= \espcondF{Y_t^{(n)}}{t}\\
  &=\expect^{\cF_{{t}}}\biggl(
        \xi+\int_{t}^T f(s,X_s,Y_s^{(n)},\widetilde{Z}_s^{(n)})\,ds
        -\int_t^T Z_s^{(n)}\,d\wien_s \\
  &\phantom{=\xi+\int_{t}^T f(s,X_s,Y_s^{(n)},\widetilde{Z}_s^{(n)})\,ds}
        -\int_{t}^{t+1/n} f(s,X_s,Y_s^{(n)},\widetilde{Z}_s^{(n)})\,ds\biggr)\\
  &=\espcondF{\xi+\int_{t+1/n}^T f(s,X_s,Y_s^{(n)},\widetilde{Z}_s^{(n)})\,ds}{t}
    -\espcondF{ {\VV}_T^{(n)}-{\VV}_t^{(n)} }{t}\\
  &=\espcondF{\xi+\int_{t+1/n}^T f(s,X_s,Y_s^{(n)},\widetilde{Z}_s^{(n)})\,ds}{t},
 \end{align*}
which proves \eqref{eq:BSDE-Ftn}.

Now, using \eqref{eq:BSDE-Ftn} and \eqref{eq:BSDE-genen}, we have
\begin{multline*}
  \espcondF{\xi+\int_0^T
    f(s,X_s,Y_s^{(n)},\widetilde{Z}_s^{(n)})\,ds}{t} \\
\begin{aligned}
 & =\espcondF{\xi+\int_{t+1/n}^T
   f(s,X_s,Y_s^{(n)},\widetilde{Z}_s^{(n)})\,ds}{t}\\
 &\phantom{f(s,X_s,Y_s^{(n)},\widetilde{Z}_s^{(n)})\,ds}
  +{\int_0^t f(s,X_s,Y_s^{(n)},\widetilde{Z}_s^{(n)})\,ds}+U_t^{(n)}
\end{aligned}
\allowdisplaybreaks\\
\begin{aligned}
  &=Y_t^{(n)}+\int_0^t
  f(s,X_s,Y_s^{(n)},\widetilde{Z}_s^{(n)})\,ds+U_t^{(n)}\\
  &={\xi+\int_{t}^T f(s,X_s,Y_s^{(n)},\widetilde{Z}_s^{(n)})\,ds
        -\int_t^T Z_s^{(n)}\,d\wien_s 
        -U_t^{(n)}}\\
   &\phantom{=\xi+}
   +\int_0^t f(s,X_s,Y_s^{(n)},\widetilde{Z}_s^{(n)})\,ds+U_t^{(n)}
\end{aligned}\\
  =\xi+\int_0^T f(s,X_s,Y_s^{(n)},\widetilde{Z}_s^{(n)})\,ds
                            -\int_t^T Z_s^{(n)}\,d\wien_s. 
\end{multline*} 
In particular,
\begin{multline*}
  \expect\CCO{\xi+\int_0^T f(s,X_s,Y_s^{(n)},\widetilde{Z}_s^{(n)})\,ds}\\
  =\xi+\int_0^T f(s,X_s,Y_s^{(n)},\widetilde{Z}_s^{(n)})\,ds
                             -\int_0^T Z_s^{(n)}\,d\wien_s. 
\end{multline*}
Thus
\begin{multline*}
  \int_0^tZ_s^{(n)}\,d\wien_s
   =\int_0^TZ_s^{(n)}\,d\wien_s-\int_t^TZ_s^{(n)}\,d\wien_s\\ 
\begin{aligned}
&=\CCO{\xi+\int_0^T f(s,X_s,Y_s^{(n)},\widetilde{Z}_s^{(n)})\,ds
      -\expect\CCO{\xi+\int_0^T f(s,X_s,Y_s^{(n)},\widetilde{Z}_s^{(n)})\,ds}}\\
&\phantom{=}-\CCO{\xi+\int_0^T f(s,X_s,Y_s^{(n)},\widetilde{Z}_s^{(n)})\,ds
      -\espcondF{\xi+\int_0^T f(s,X_s,Y_s^{(n)},\widetilde{Z}_s^{(n)})\,ds}{t}}\\
&=\espcondF{ \xi+\int_0^T f(s,X_s,Y_s^{(n)},\widetilde{Z}_s^{(n)})\,ds}{t}
-\expect\CCO{\xi+\int_0^T f(s,X_s,Y_s^{(n)},\widetilde{Z}_s^{(n)})\,ds}
\end{aligned}
\end{multline*}
which proves \eqref{eq:BSDE-Zn}.
\finpr

\subsection{Boundedness and continuity results}
\label{subsect:boundedness}
In this part, we show some results that will be useful to
 prove the relative compactness in distribution of the sequence 
$\bigl(Y^{(n)},
\int_{.}^T f(s,X_s,Y_s^{(n)},\widetilde{Z}_s^{(n)})\,ds,\allowbreak
\int_{.}^T{Z}_s^{(n)}\,d\wien_s,\allowbreak
Z^{(n)}\bigr)$ 
in some properly chosen state space.

\begin{lem}\label{lem:Zbounded}
Let
\begin{multline*}
 \widetilde{Y}_t^{(n)} 
      = {Y}_t^{(n)} + U_t^{(n)}
      = {\xi}+\int_t^{T} f(s,X_s,Y_s^{(n)},\widetilde{Z}_s^{(n)})\,ds
                             -\int_t^{T} Z_s^{(n)}\,d\wien_s.
\end{multline*}
There exist constants $\CA,\CB>0$ such that, for all $t$ such that 
$0\leq t\leq T$,  
\begin{equation}
  \label{eq:Zbounded}
   \expect\int_t^{T}\norm{Z_s^{(n)}}^2ds
     \leq \CA\expect\int_t^{T}\norm {\widetilde Y_s^{(n)}}^2ds + \CB.
\end{equation}
\end{lem}
\proof 
Using Proposition \ref{prop:YnZn_inegrables}, we have, for each $n\geq 1$, 
\begin{multline}\label{eq:supYtilde}
\expect \CCO{\sup_{t\in[0,T]}\norm{\widetilde{Y}_t^{(n)}}^2}\\
\leq 3\expect\norm{\xi}^2+3\consf^2\expect\int_0^T
\bigl(1+\norm{Z_s^{(n)}}\bigr)^2 ds
+3\expect \CCO{\sup_{t\in[0,T]}\norm{\int_t^{T} Z_s^{(n)}\,d\wien_s}^2}
<+\infty.
\end{multline} 
Applying It\^o's formula to the semi-martingale $\norm
{\widetilde{Y}_t^{(n)}}^2$, taking expectation of both sides and using
the fact that 
$t'\mapsto \int_t^{t'} \scal{\widetilde{Y}_s^{(n)}}{
  Z_s^{(n)}\,d\wien_s}$ is a martingale 
(thanks to \eqref{eq:supYtilde} and Proposition \ref{prop:YnZn_inegrables}), 
we get 
$$\expect\norm{\widetilde{Y}_t^{(n)}}^2
= \expect \norm{\xi}^2
     +2\expect \int_t^{T}
           {\widetilde{Y}_s^{(n)}}.{f(s,X_s,Y_s^{(n)},\widetilde{Z}_s^{(n)})}ds
- \expect \int_t^{T} \norm{Z_s^{(n)}}^2ds.$$
Thus
\begin{align*}
 \expect \int_t^{T} \norm{Z_s^{(n)}}^2\,ds
\leq \expect \norm {\xi}^2
     +2\expect \int_t^{T}                     
        \norm{\widetilde{Y}_s^{(n)}}. 
            \norm{f(s,X_s,Y_s^{(n)},\widetilde{Z}_s^{(n)})}\,ds.
\end{align*}
From \Ha, this entails
\begin{align*}
  \expect \int_t^{T} \norm{Z_s^{(n)}}^2\, ds
 &\leq \expect \norm {\xi}^2 
  + 2 \consf \expect \int_t^{T}  \norm{\widetilde{Y}_s^{(n)}} 
                                   (1+\norm{\widetilde{Z}_s^{(n)}})\,ds.
  \end{align*}
Using that, for $a\geq0, b\geq0$, and $\lambda\not=0$, 
we have $2ab\leq a^2\lambda ^2 + b^2/\lambda ^2$, we get 
\begin{multline*}
2\expect\int_t^T \norm{\widetilde Y_s^{(n)}}(1+\norm{\widetilde{Z}_s^{(n)}})\,ds\\
\begin{aligned}
&\leq \lambda^2 \expect\int_t^T\norm{\widetilde Y_s^{(n)}}^2ds
      +2(T-t)/\lambda^2
      +2/\lambda ^2\expect\int_t^T\norm{\widetilde{Z}_s^{(n)}}^2ds\\
&\leq \lambda^2 \expect\int_t^T\norm{\widetilde Y_s^{(n)}}^2ds
      +2(T-t)/\lambda^2
      +2/\lambda ^2\expect\int_t^T\norm{ {Z}_s^{(n)}}^2ds.
\end{aligned}
\end{multline*}
Thus, taking ${\lambda^2}>2 \consf$,
\begin{align*}
(1-2 \consf/\lambda^2)\expect \int_t^{T} \norm{Z_t^{(n)}}^2\,ds
    &\leq  \expect \norm {\xi}^2 
  + \consf \CCO{
       2 T/\lambda^2 
           + \lambda^2\expect \int_t^{T} \norm{\widetilde{Y}_t^{(n)}}^2\,ds}
\end{align*}
which yields \eqref{eq:Zbounded}.
\finpr

\begin{prop}\label{prop:YZbounded}
  Let $\widetilde{Y}_t^{(n)} = {Y}_t^{(n)} + U_t^{(n)}$ be as in Lemma
  \ref{lem:Zbounded}. The families 
  $( \widetilde Y_t^{(n)})_{0\leq t\leq T,\, n\geq 1}$, $(
  Y_t^{(n)})_{0\leq t\leq T,\, n\geq 1}$ and 
$(U_t^{(n)})_{0\leq t\leq T,\, n\geq 1}$ are bounded in 
$\ellp{2}_{\espY}(\esprob)$. 
\end{prop}
\proof
We have
\begin{align*}
\widetilde Y_t^{(n)}
                    &= {Y}_t^{(n)} + U_t^{(n)}\\
                    &=\espcondF{
              \xi+\int_{t+1/n}^T f(s,X_s,Y_s^{(n)},\widetilde{Z}_s^{(n)})\                                     }{t}
                         + \espcondF{ 
               \int_t^{t+1/n} f(s,X_s,Y_s^{(n)},\widetilde{Z}_s^{(n)})\,ds 
                                   }{t}\\
                   &= \espcondF{ \xi 
              +\int_{t}^T f(s,X_s,Y_s^{(n)},\widetilde{Z}_s^{(n)})\,ds
                                  }{t}.
\end{align*}
We deduce the following inequalities, where $\consp$ denotes some
constant which is not necessarily the same at each line 
but does not depend on $n$:
\begin{align*}
  \expect\norm{\widetilde Y_t^{(n)}}^2
   &=  \expect \norm{
                \espcondF{
                   \xi+\int_{t}^T f(s,X_s,Y_s^{(n)},\widetilde{Z}_s^{(n)})\,ds
                          }{t}
                                }^2\\
 &\leq \consp \expect \CCO{
           {\norm{\xi}^2}
           +\int_t^T(1+\norm{Z_s^{(n)}}^2)\,ds 
                    }\\
 &\leq \consp \CCO{ 1 +\int_t^T\expect\norm{\widetilde Y_s^{(n)}}^2 ds  }.
                                          \label{eq:pre-Gronwall}
\end{align*}
The last inequality is a consequence of Lemma \ref{lem:Zbounded}. 
Let $g(t)=\expect\norm{\widetilde Y_{T-t}^{(n)}}^2$. 
The preceding inequalities yield 
\begin{align*}
  g(t)&\leq \consp\CCO{ 1+\int_0^t g(s)\,ds }.
\end{align*}
Thus, by Gronwall's Lemma, 
\begin{align*}
  g(t)&\leq \consp \CCO{1+\consp \int_0^t e^{\consp(t-s)}\,ds}
       \leq \consp \CCO{1+\consp \int_0^T e^{\consp(T-s)}\,ds}
\end{align*}
which proves that $(\widetilde Y_t^{(n)})_{0\leq t\leq T,\, n\geq 1}$ 
  is bounded in $\ellp{2}_{\espY}(\esprob)$.\\ 
Now, we have, using again Lemma \ref{lem:Zbounded}, 
\begin{align*}
 \expect\CCO{\norm{ Y_t^{(n)}}^2}
   &=  \expect \norm{
                  \espcondF{
             \xi+\int_{t+1/n}^T f(s,X_s,Y_s^{(n)},\widetilde{Z}_s^{(n)})\,ds
                          }{t}
                                }^2\\
 &\leq \consp \expect \CCO{
           {\norm{\xi}^2}
           +\int_{t+1/n}^T(1+\norm{Z_s^{(n)}}^2)\,ds 
                    }\\
 &\leq \consp \CCO{ 1 +\int_{t+1/n}^T\expect\norm{\widetilde Y_s^{(n)}}^2\,ds  }
\end{align*}
which proves that $(Y_t^{(n)})_{0\leq t\leq T,\, n\geq 1}$  is bounded
in $\ellp{2}_{\espY}(\esprob)$. 

The boundedness in $\ellp{2}_{\espY}(\esprob)$ of  
$(U_t^{(n)})_{0\leq t\leq T,\, n\geq 1}$ 
follows immediately from that
of $(\widetilde{Y}_t^{(n)})_{0\leq t\leq T,\, n\geq 1}$ 
and $({Y}_t^{(n)})_{0\leq t\leq T,\, n\geq 1}$.
\finpr

\begin{cor}\label{cor:Zbounded} The sequences 
$(Z^{(n)})_{n\geq 1}$ and $(\widetilde{Z}^{(n)})_{n\geq 1}$
  are bounded in $\ellp{2}_{\espZ}(\esprob\times[0,T])$, and we have 
\begin{equation}\label{eq:supVbounded}
\sup_{n\geq 1}\expect\CCO{
\sup_{0\leq t\leq T}\norm{\int_t^T Z_s^{(n)}\,d\wien_s
                   }^2}
<+\infty.
\end{equation}
\end{cor}
\proof
The boundedness of $(Z^{(n)})_{n\geq 1}$ and
$(\widetilde{Z}^{(n)})_{n\geq 1}$ 
is a direct consequence of 
Lemma \ref{lem:Zbounded} and Proposition \ref{prop:YZbounded}.
Then \eqref{eq:supVbounded} follows by It\^o's isometry, Doob's
inequality, and the fact that
$$\norm{\int_t^T Z_s^{(n)}\,d\wien_s}
\leq\norm{\int_0^T Z_s^{(n)}\,d\wien_s}+\norm{\int_0^t Z_s^{(n)}\,d\wien_s}. $$ 
\finpr

\begin{lem}\label{lem:Utn}
Let $1\leq \pq<2$.
We have 
\begin{equation}
\lim_{n\rightarrow\infty}
\expect\biggl(\sup_{0\leq t\leq T}\norm{U_t^{(n)}}^\pq\biggr)=0.
                  \label{eq:Utn-0} 
\end{equation}
\end{lem}
\proof
For each $n$, we can find 
an $\cF_T$-measurable time $\tau_n$ such that 
\begin{equation*}
{\sup_{0\leq t\leq T}
  {\int_t^{t+1/n}\CCO{1+\norm{Z_s^{(n)}}}^\pq ds}}
={\int_{\tau_n}^{\tau_n+1/n}\CCO{1+\norm{Z_s^{(n)}}}^\pq ds}.
\end{equation*}

Let $\majorr{2}=\sup_{n}\expect\int_0^T\CCO{1+\norm{Z_s^{(n)}}}^2ds$. 
By 
Corollary \ref{cor:Zbounded}, 
we have $\majorr{2}<+\infty$. Let $\pq'$ such that $\pq<\pq'<2$. 
Using the growth condition {\Ha} and Doob's inequality applied to the
martingale 
$\espcondF{\int_{\tau_n}^{\tau_n+1/n}\CCO{1+\norm{Z_s^{(n)}}}^\pq ds}{t}$, 
we get 
\begin{align*}
\expect\CCO{\sup_{0\leq t\leq T}\norm{U_t^{(n)}}^\pq}
&\leq \consf^\pq\expect\CCO{\sup_{0\leq t\leq T}
  \espcondF{\int_t^{t+1/n}\CCO{1+\norm{Z_s^{(n)}}}^\pq ds}{t}
                     }\\
&\leq \consf^\pq\expect \CCO{\sup_{0\leq t\leq T}
  \espcondF{\int_{\tau_n}^{\tau_n+1/n}\CCO{1+\norm{Z_s^{(n)}}}^\pq ds}{t}
                     }\\
&\leq \consf^\pq\CCO{\expect \CCO{\sup_{0\leq t\leq T}
  \espcondF{{\int_{\tau_n}^{\tau_n+1/n}\CCO{1+\norm{Z_s^{(n)}}}^\pq ds}}{t}^{\pq'/\pq}
                     }}^{\pq/\pq'}\\
&\leq \frac{\pq'}{\pq'-\pq}\,\consf^\pq\CCO{\expect{
  \CCO{\int_{\tau_n}^{\tau_n+1/n}\CCO{1+\norm{Z_s^{(n)}}}^{\pq'} ds} }
                     }^{\pq/\pq'}\\
&\leq\frac{\pq'}{\pq'-\pq}\, \consf^\pq
  \CCO{\frac{1}{n}}^{(2-\pq')/2}
  \CCO{\CCO{
                 \expect\int_0^{T}\CCO{ 1+\norm{Z_s^{(n)}} }^2
                 ds}^{\pq'/2}
          }^{\pq/\pq'}
                     \\
&=\frac{\pq'}{\pq'-\pq}\, {\consf^\pq{
  \CCO{\frac{1}{n}}^{(2-\pq')/2}\CCO{
                 \expect\int_0^{T}\CCO{1+\norm{Z_s^{(n)}}}^2
                 ds}^{\pq/2}}
                     }\\
&\leq\frac{\pq'}{\pq'-\pq}\,
\consf^\pq\,{\majorr{2}^{\pq/2}\,\CCO{\frac{1}{n}}^{(2-\pq')/2}
               },
\end{align*}
which proves \eqref{eq:Utn-0}. 
\finpr

\begin{lem}\label{lem:supYbounded}
  We have 
 \begin{equation*}
\sup_{n\geq 1}\,\expect\CCO{\sup_{0\leq t\leq T}\norm{Y_t^{(n)}}^2}
<+\infty.
\end{equation*}
\end{lem}
\proof
Using \eqref{eq:BSDE-genen}, we get
\begin{align*}
          \sup_{0\leq t\leq T}\norm{Y_t^{(n)}}^2
  &\leq A_n+B_n+C_n
\end{align*}
where
\begin{align*}
A_n&=3\sup_{0\leq t\leq T}
\norm{\xi+\int_{t+1/n}^Tf(s,X_s,Y_s^{(n)},\widetilde{Z}_s^{(n)})\,ds}^2,\\ 
B_n&=3\sup_{0\leq t\leq T} \norm{\int_t^T Z_s^{(n)}\,d\wien_s}^2, \\ 
  C_n&=3\sup_{0\leq t\leq T} \norm{U_t^{(n)} }^2. 
\end{align*}
By 
Corollary \ref{cor:Zbounded}, 
$(Z^{(n)})_{n\geq 1}$ 
  is bounded in $\ellp{2}_{\espZ}(\esprob\times[0,T])$, thus 
  using
  the growth condition \Ha, we get 
\begin{align*}
\sup_n \expect\CCO{\sup_{0\leq t\leq
    T}\biggl(\norm{\xi}^2
   +\consf^2\int_{t+1/n}^T\CCO{1+\norm{Z_s^{(n)}}}^2\,ds\biggr)
  }<+\infty 
\end{align*}
which entails $\sup_n \expect(A_n)<+\infty$. 
On the other hand, 
$ \VV_t^{(n)}:=\int_0^t Z_s^{(n)}\,d\wien_s$  is a martingale, 
so, using again 
Corollary \ref{cor:Zbounded},
\begin{align*}
 \sup _n \expect\CCO {B_n}
    & \leq  \consp\sup_n \expect\norm{\VV_T^{(n)}}^2
<+\infty.
\end{align*}
Finally from 
\Ha\ 
and the boundedness of $(Z^{(n)})_{n\geq 1}$ 
in $\ellp{2}_{\espZ}(\esprob\times[0,T])$ (see Corollary \ref{cor:Zbounded}),
we have 
\begin{align*}
\sup_n \expect \CCO {C_n}
&\leq 3\sup_{0\leq t\leq T}\consf^2
    \expect\CCO{{\int_0^T \CCO{1+\norm{Z_s^{(n)}}}^2\,ds}}
<+\infty.
\end{align*} 
\finpr

\subsection{Compactness results}\label{subsect:compactness}

\begin{lem}\label{lem:ftight}
  The sequence 
  $(\int_{.}^T f(s,X_s,Y_s^{(n)},\widetilde{Z}_s^{(n)})\,ds
)_{n\geq 1}$ 
  is tight in $\trajY$. 
\end{lem}
\proof
Let us denote $\FYZ^{(n)}=\int_{.}^T
f(s,X_s,Y_s^{(n)},\widetilde{Z}_s^{(n)})\,ds$. 
By a criterion of Aldous 
\cite{aldous78stopping,jacod85limite}, we only need to prove that 
 $$
 \forall \epsilon>0, \ \exists R >0, \forall n\geq 1,\ 
 \prob\CCO{\sup_{0\leq t\leq T}\norm{\FYZ_t^{(n)}}\geq R}\leq \epsilon
 \leqno{(A)}$$
$$ 
\forall \epsilon>0, \forall \eta>0, \ \exists \delta>0: \forall n\geq 1, 
\sup_{\substack{\sigma,\tau\in\TA\\
                                0\leq\abs{\tau-\sigma}\leq\delta}}
  \prob\CCO{\norm{\FYZ_\tau^{(n)} -\FYZ_ \sigma^{(n)}}\geq \eta }\leq \epsilon
     \leqno{(B)}$$
where $\TA$ denotes the set of stopping times with values in $[0,T]$. 
We are going to prove the slightly stronger properties 
\begin{gather}
\sup_{n\geq 1}\,\expect\CCO{\sup_{0\leq t\leq T}\norm{\FYZ_t^{(n)}}}
<+\infty,\label{eq:Aprime}\\
\label{eq:Bprime}
\forall \epsilon>0, 
\ \exists \delta>0: 
\sup_{n\geq 1}\,\sup_{\substack{\sigma,\tau\in\TA\\
               \abs{\tau-\sigma}\leq\delta}}
              \expect{\norm{\FYZ_\tau^{(n)} -\FYZ_ \sigma^{(n)}}}
<+\epsilon.
\end{gather}
As 
$$\FYZ_t^{(n)}=Y_t^{(n)}-\xi
+\int_t^T Z_s^{(n)}\,d\wien_s 
        +U_t^{(n)},$$
we can, for example, deduce \eqref{eq:Aprime} from 
Corollary \ref{cor:Zbounded}, 
Lemma \ref{lem:Utn}, and
Lemma \ref{lem:supYbounded}.

Now, let $\sigma,\tau\in\TA$, with $\abs{\tau-\sigma}\leq\delta$. 
We have 
$\expect{\norm{\FYZ_\tau^{(n)} -\FYZ_ \sigma^{(n)}}}=
\expect{\norm{\FYZ_{\sigma\vee\tau}^{(n)} -\FYZ_{\sigma\wedge\tau}^{(n)}}}$. 
Thus  we can assume without loss of generality
that $\sigma\leq \tau$. 
Then
\begin{align*}
\expect{\norm{\FYZ_\tau^{(n)} -\FYZ_ \sigma^{(n)}}}
&=
\expect\norm{
\int_{\sigma}^{\tau} f(s,X_s,Y_s^{(n)},\widetilde{Z}_s^{(n)})\,ds}\\ 
&\leq \CCO{\expect(\tau-\sigma)}^{1/2}  \CCO{\expect 
\int_0^T \norm{ f(s,X_s,Y_s^{(n)},\widetilde{Z}_s^{(n)})}^2\,ds
                                                     }^{1/2}\\
&\leq\delta^{1/2}\consf\CCO{\expect{\int_0^T}(1+\norm{Z_s^{(n)}})^2\,ds}^{1/2}
\end{align*}
and \eqref{eq:Bprime} follows from Corollary  \ref{cor:Zbounded}.
\finpr

\paragraph{The topology S and Condition UT}
In order to prove the tightness of  $(Y^{(n)})_{n\geq 1}$, 
we will use Meyer-Zheng criterion 
\cite{Meyer-Zheng84tightness}
and Jakubowski's topology S \cite{jakubowski97skor} on the space 
$\trajYD:=\trajYD_\espY[0,T]$. 
First, we need some
definitions.

Let $\trajYV\subset\trajYD$ be the subspace  
of elements of $\trajYD$ which have finite variation. 
The topology S on $\trajYD$ 
is defined by its convergent sequences: A sequence 
$(x_n)$ in $\trajYD$ converges for S to a limit
$x\in\trajYD$ if, from any subsequence of $(x_n)$, one can extract a
further subsequence $(x'_n)$ such that, for every $\epsilon>0$, there
exist a sequence $(v_{n,\epsilon})$ of elements of $\trajYV$
and $v_{\epsilon}\in\trajYV$ (depending on the subsequence $(x'_n)$) 
such that 
\begin{enumerate}[(i)]
\item $\sup_n  
\sup_{t\in[0,T]}\norm{x'_n(t)-v_{n,\epsilon}(t)}\leq\epsilon$ and 
$\sup_{t\in[0,T]}\norm{x(t)-v_{\epsilon}(t)}\leq\epsilon$,
\item $\lim_{n\rightarrow\infty}
\int_0^T f(t)\,dv_{n,\epsilon}(t)
=\int_0^T f(t)\,dv_{\epsilon}(t)$
for every continuous function $f$ defined on $[0,T]$. 
\end{enumerate}
We denote $\trajYDs$ the space $\trajYD$ endowed with S. 
The topology S is coarser than Skorokhod's topology $J_1$, which is
Polish, thus S is Lusin (see \cite{schwartz73book} on 
properties of Lusin spaces). 
In particular, by \cite[Corollary 2 page 101]{schwartz73book}, 
S has the same
Borel sets as $J_1$, 
thus the Borel subsets of S are 
generated by the projection
mappings $\proj{t} :\, x\mapsto x(t)$ for $t\in[0,T]$). 
Furthermore, S is finer than 
the Meyer-Zheng topology \cite{Meyer-Zheng84tightness}, 
which is the topology on $\trajYV_\espY$ induced by $\ellp{0}_\espY([0,T],dt)$.
In particular, 
S is (separably) submetrizable, that is, 
there exists a (separable) metrizable topology which
is coarser than S. Equivalently, one can find a countable set of
S-continuous real-valued functions  which separate the points
of $\trajYD$. 
This implies that S is Hausdorff and that 
the compact subsets of
$\trajYDs$ are metrizable. 

Another important feature of $S$ is that the addition $(x,y)\mapsto
x+y$ is S-sequentially continuous on $\trajYDs\times\trajYDs$.

A criterion of tightness on $\trajYDs$ is the
so-called condition UT (see \cite[Theorem 4.2]{jakubowski97skor}):
Let $\mathcal{H}$ denote the set of elementary
real valued predictable processes bounded by 1, i.e.~processes of the
form 
 $$H_t=\un{[t_0,t_1]}(t)H_0+\un{]t_1,t_2]}(t)H_1+\dots
+\un{]t_{n-1},t_{n}]}(t)H_{t_{n-1}}$$
where $0=t_0\leq\dots\leq t_n\leq T$ and each $H_i$ is bounded by 1
and $\cF_{t_i}$-measurable. 
Let $(\procgenne^\alpha)_{\alpha\in A}$ be a family 
of $\trajYD$-valued processes. 
We say that $(\procgenne^\alpha)$
{\em satisfies Condition UT}\, if 
the family of all stochastic
integrals $\int H\,d\procgenne^\alpha$, 
where $\alpha\in A$ and $H\in\mathcal{H}$,
is uniformly tight. 
Condition UT was 
considered for the first time by Stricker \cite{stricker85compacite},
to prove compactness in the Meyer-Zheng topology.  
Discussions on this condition can be found in 
\cite{jakubowskimeminpages89conv,Memin-Slominski91UT}. 

We now consider a stronger condition, 
proposed by Meyer and Zheng \cite{Meyer-Zheng84tightness}: 
Let  $\procgene$ be an adapted process defined on the time interval
$[0,T]$, with values in $\espY$. 
For any finite partition 
$\pi=(t_0,\dots,t_n)$ of $[0,T]$, let us denote
$$\cvar_\pi\CCO{\procgene}
=\expect{\norm{\procgene_T}+\sum_{i=0}^{n-1}
\norm{ \espcondF{ \procgene_{t_{i+1}}-\procgene_{t_{i}}  }{{t_i}} }
},$$
and define the {\em conditional variation} $\cvar\CCO{\procgene}$ of
$\procgene$ by 
$$\cvar\CCO{\procgene}=\sup_\pi\cvar_\pi\CCO{\procgene}.$$
By \cite[Th\'eor\`eme 3]{stricker85compacite}, 
if a family $(\procgenne^\alpha)$ of adapted
$\trajYD$-valued processes 
satisfies 
$$\sup_{\alpha}\cvar\CCO{\procgenne^\alpha}<\infty,$$
then  Condition UT holds for $(\procgenne^\alpha)$. 

An adapted stochastic process $\procgenne$ such that
$\cvar\CCO{\procgene}<\infty$ is called a {\em quasimartingale}. 
Let us mention that, if the quasimartingale $\procgenne$ is
right-continuous in probability, then it has a \cadlag\ adapted
version 
(assuming the right-continuity of $(\cF_t)$), 
see \cite[Theorem 4.1]{brooks-dinculeanu87quasimartingales}.

\begin{prop}\label{lem:Ytight}
  The sequences 
 $(Y^{(n)})_{n\geq 1}$ and $(\int_{.}^TZ_s^{(n)}\,d\wien_s)_{n\geq 1}$  
  are tight sequences of $\trajYD$-valued random variables, for the
  topology S. 
\end{prop}
\proof
First, we need to check that, for each integer $n\geq 1$, the process
\begin{equation}\label{eq:mapF}
Y^{(n)}_t=\espcondF{\xi
      +\int_{t+1/n}^T f(s,X_s,Y_s^{(n)},\widetilde{Z}_s^{(n)})\,ds}{t}
\end{equation}
has a $\trajYD$-valued version. 
Let us prove that it is continuous in
$\ellp{1}$ and a quasimartingale.  
As $(\cF_t)$ is a Brownian filtration, the
martingale 
$$t\mapsto
\espcondF{\xi
      +\int_{r+1/n}^T f(s,X_s,Y_s^{(n)},\widetilde{Z}_s^{(n)})\,ds}{t}
$$
has a continuous version 
for each fixed $r\in[0,T-1/n]$, thus it is continuous in
$\ellp{1}$, i.e.~the mapping
\begin{equation}\label{eq:mapping}
\left\{
\begin{array}{lcl}
[0,T]\times[0,T-1/n]&\rightarrow&\ellp{1}\\
(t,r)&\mapsto&\espcondF{\xi
      +\int_{r+1/n}^T f(s,X_s,Y_s^{(n)},\widetilde{Z}_s^{(n)})\,ds}{t}
\end{array}
\right.
\end{equation} is continuous in the variable $t$. 
On the other hand, we have, 
for fixed $t\in[0,T]$ and for $0\leq r_1\leq r_2\leq T-1/n$ such that
$r_2-r_1\leq 1/n$,
\begin{multline*}
\norm{\espcondF{
    \int_{r_1+1/n}^{r_2+1/n}f(s,X_s,Y_s^{(n)},\widetilde{Z}_s^{(n)})\,ds}{t}}\\
\begin{aligned}
&\leq \consf\espcondF{
     \int_{r_1+1/n}^{r_2+1/n}(1+\norm{\widetilde{Z}_s^{(n)}})\,ds}{t}\\
&\leq (r_2-r_1)^{1/2}\consf\CCO{\espcondF{
     \int_{0}^{T}(1+\norm{\widetilde{Z}_s^{(n)}})^2\,ds}{t}}^{1/2}.
\end{aligned}
\end{multline*}
Therefore, by Corollary \ref{cor:Zbounded}, the mapping 
\eqref{eq:mapping}
is continuous in $r$ uniformly with respect to $t$, 
thus it is jointly continuous, which proves 
the continuity in $\ellp{1}$ of the process \eqref{eq:mapF} for each
$n\geq 1$.

Now, we have, for any subdivision $\pi=(t_0,\dots,t_m)$ of $[0,T]$,  
\begin{align*}
\sup_{n}\,\cvar_\pi(Y^{(n)})
=&\sup_{n}\,\expect\CCO{\norm{\xi}+\sum_{i=0}^{m-1}
\norm{ \espcondF{  Y^{(n)}_{t_{i+1}}-Y^{(n)}_{t_{i}}  }{{t_i}} }
}\\
\leq& \sup_{n}\,\expect\CCO{\norm{\xi}+\sum_{i=0}^{m-1}
\norm{ 
\int^{t_{i+1}+1/n}_{t_{i}+1/n}
f(s,X_s,Y_s^{(n)},\widetilde{Z}_s^{(n)})\,ds 
     }
}\\
\leq& \sup_{n}\,\expect\CCO{\norm{\xi}+
\int_0^T \norm{f(s,X_s,Y_s^{(n)},\widetilde{Z}_s^{(n)})}\,ds
}.\\
\intertext{This estimation does not depend on $\pi$, thus, 
using Corollary \ref{cor:Zbounded},}
\sup_{n}\,\cvar(Y^{(n)})
\leq& \expect\CCO{\norm{\xi}}+\sup_{n}\,\expect\CCO{{
\int_0^T \consf(1+\norm{\widetilde{Z}_s^{(n)}})\,ds
}}\\
\leq& \expect\CCO{\norm{\xi}}+\sup_{n}\, \consf\CCO{T+T^{1/2}
\CCO{\expect\CCO{\int_0^T\norm{\widetilde{Z}_s^{(n)}}^2)\,ds }}^{1/2}}\\
<&+\infty.
\end{align*}
This proves that each $Y^{(n)}$
is a quasimartingale, 
and that the sequence $(Y^{(n)})_{n\geq 1}$ satisfies
Condition UT. 
Furthermore, for each $n\geq 1$, as $Y^{(n)}$ is right-continuous in
$\ellp{1}$, it has 
a \cadlag\ version, 
thanks to \cite[Theorem 4.1]{brooks-dinculeanu87quasimartingales}.
Thus, by \cite[Theorem 4.2]{jakubowski97skor},
the sequence $(Y^{(n)})_{n\geq 1}$ is tight in $\trajYDs$.

Similarly,
\begin{multline*}
\sup_{n}\,\cvar_\pi\biggl(\int_{0}^.Z_s^{(n)}\,d\wien_s\biggr)\\
\begin{aligned}
&=
\sup_{n}\,\expect
\biggl(\norm{\int_{0}^{T}  Z_s^{(n)}\,d\wien_s} 
+\sum_{i=0}^{m-1}
\norm{ \espcondF{ \int_{t_{i}}^{t_{i+1}}  Z_s^{(n)}\,d\wien_s  }{{t_i}} }
\biggr)\\
&=\sup_{n}\,\expect
\norm{\int_{0}^{T}  Z_s^{(n)}\,d\wien_s},
\end{aligned}
\end{multline*} 
{thus}
\begin{align*}
\sup_{n}\,\cvar\biggl(\int_{0}^.Z_s^{(n)}\,d\wien_s\biggr)
=&\sup_{n}\,\expect
\norm{\int_{0}^{T}  Z_s^{(n)}\,d\wien_s} 
<+\infty.
\end{align*} 
Thus $(\int_{0}^{.}Z_s^{(n)}\,d\wien_s)_{n\geq 1}$ satisfies Condition
UT. 
Again 
by \cite[Theorem 4.2]{jakubowski97skor}, this proves that 
$(\int_{0}^{.}Z_s^{(n)}\,d\wien_s)_{n\geq 1}$ is tight in $\trajYDs$. 
Finaly, it is straightforward to check that the mapping 
$$\left\{
\begin{array}{lcl}
\trajYD&\rightarrow&\trajYD\\
u&\mapsto&u(T)-u
\end{array}
\right.$$ 
is sequentially continuous for the topology S, thus 
$(\int_{.}^TZ_s^{(n)}\,d\wien_s)_{n\geq 1}$ is tight in $\trajYDs$. 
\finpr

\subsection{Construction of a weak limit process}
\label{subsect:weak_solution}
This part of 
the construction of a weak solution 
follows the same lines as in
\cite{jakubowski-kamenski-prf05stochinc}, with some complications due
to the processes $Z^{(n)}$.

\paragraph{Young measures}
Let us recall the definition and main properties of Young measures, 
see \cite{valadier94course, balder00lectures} for 
introductions to the topic, 
and \cite{cc-prf-valadier04book} for the setting of nonnecessarily
regular topological spaces, which we need here.  
Let $\espE$ be a 
Suslin topological space 
(i.e.~$\espE$ is a Hausdorff topological space and 
there exists a Polish space $\esppolish$ and a continuous surjective mapping
from $\esppolish$ onto $\espE$, 
see \cite{schwartz73book} for the properties of Suslin spaces, 
or \cite[Chapter 1]{cc-prf-valadier04book} for a survey without proofs). 
Let $\bor{\espE}$ be the Borel
$\sigma$-algebra of $\espE$.  
A {\em Young measure $\mu$ with basis $\prob$ on $\espE$} is a
probability measure on $\esprob\times \espE$,   
such that for any set $A\in \cF$, $\mu(A\times \espE)=\prob(A)$. The
space of Young measures with basis $\prob$ is denoted by
$\youngs(\esprob,\tribu, \prob; \espE)$. 
It is very useful to describe a Young measure $\mu$ by its {\em
  disintegration} $(\mu_\omega)$ with respect to $\prob$ 
(see Definition \ref{def:disintegration}). 
 The space $\ellp{0}(\esprob,\tribu, \prob;\espE)$ of measurable functions from
 $\esprob$ to $\espE$ is embedded in $\youngs(\esprob,\tribu, \prob;
 \espE)$ in the following way: we identify 
every $u \in \ellp{0}(\esprob,\tribu, \prob;\espE) $ 
with the Young measure $\delta_{u(\omega)}\otimes \,d\prob(\omega)$, 
where $\delta_{u(\omega)}$ denotes the Dirac mass at ${u(\omega)}$. 
In other words, $u$ is identified with 
the unique Young measure $\mu$ whose support is the graph of $u$. 
The set $\youngs(\esprob,\tribu, \prob; \espE)$ is endowed with a
topology defined as follows: 
 A generalized sequence\footnote{see \cite{kelley55book} on
   generalized sequences, also called nets, 
 however we do not need them in the sequel, because we
 use sequential compactness results. Note also that, 
when $\espE$ is metrizable, the space 
$\youngs(\esprob,\tribu, \prob; \espE)$
is metrizable too, 
and we can characterize its topology using convergent sequences 
instead of convergent generalized sequences.}  
$(\mu^\alpha)$ 
of Young measures converges to a Young measure
 $\mu$ if, for each bounded measurable 
$\Phi :\esprob\times \espE\rightarrow\R$ such that 
$\Phi(\omega, .)$ is continuous for all $\omega\in\esprob$, 
the generalized sequence $(\mu^\alpha(\Phi))$ converges to
$\mu(\Phi)$. 
In this case, we say that $(\mu^\alpha)$ {\em converges stably}, 
or {\em $\tribu$-stably}, 
to $\mu$ (this terminology stems from R\'enyi \cite{renyi63stable}).

Note that the restriction of the topology of stable convergence to
$\ellp{0}(\esprob;\espE)$ is the topology of convergence in
probability, see \cite{valadier94course,cc-prf-valadier04book}. 

We say that a subset $\mathcal{K}$ of $\youngs(\esprob,\tribu, \prob;
\espE)$ is {\em tight} if, for each $\epsilon>0$, there exists a
compact subset $K$ of $\espE$ such that 
$\inf_{\mu\in\mathcal{K}}\mu(\esprob\times K)\geq 1-\epsilon$. In the
case when $\mathcal{K}\subset\ellp{0}(\esprob,\tribu, \prob ;\espE)$, 
this is the usual tightness notion for random variables. 
By \cite[Theorem 4.3.5]{cc-prf-valadier04book}, if the compact subsets
of $\espE$ are metrizable, and if $\mathcal{K}$ is tight, then
$\mathcal{K}$ is relatively compact and relatively sequentially
compact in $\youngs(\esprob,\tribu, \prob;\espE)$. 
The converse is true if $\espE$ has the Prohorov property. 

We will need a result on convergence of Young measures with respect to
sequentially continuous integrands:
\begin{lem}\label{rem:seq-continu}%
Assume that $\espE$ is a Suslin submetrizable topological
space. 
Let $(\mu^n)$ be a tight sequence in $\youngs(\esprob,\tribu, \prob; \espE)$
which stably converges to a Young measure $\mu$. 
Let $f :\,\esprob\times\espE \rightarrow\R$ be a bounded measurable
function such that $f(\omega,.)$ is {\em sequentially} continuous for
each $\omega\in\esprob$. Then 
$\lim_n\mu^n(f)=\mu(f)$. 
\end{lem}
\proof
By Balder's extension of Koml\'os Theorem for Young measures 
\cite{balder89Prohorov,balder90new} which is
valid for Hausdorff spaces with metrizable compact subsets 
\cite[Lemma 4.5.4]{cc-prf-valadier04book},  
we can extract from every subsequence of $(\mu^n)$ a further 
subsequence (which we still denote by $(\mu^n)$ for simplicity of
notations), which  
K--converges to $\mu$, that is, for each subsequence 
$({\nu}^{n})$ of 
$(\mu^n)$, 
we have 
\begin{equation}\label{eq:komlos}
\lim_n\frac{1}{n}\sum_{k=1}^{n}{\nu}^n_\omega
=\mu_\omega\ \text{a.e.}
\end{equation}
where the limit is taken in the narrow convergence, 
i.e.~$\lim_n\frac{1}{n}\sum_{k=1}^{n}{\nu}^n_\omega(g)=\mu_\omega(g)$
for every bounded continuous $g :\, \espE\rightarrow\R$.   
Let us denote $\lambda^n=\frac{1}{n}\sum_{k=1}^{n}{\nu}^n$, and let us
prove that 
\begin{equation}\label{eq:komlos-f}
\lim_n \int_\espE f(\omega,x)\,d\lambda^n_\omega(x)
=\int_\espE f(\omega,x)\,d\mu_\omega(x)\text{ a.e.}
\end{equation}
Let $\omega$ be in the almost sure set on which the convergence in 
\eqref{eq:komlos} holds. 
As $\espE$ admits a coarser separable metrizable topology, we can
apply Jakubowski's extension of Skorokhod's representation theorem
\cite{jakubowski97almo}: 
for every subsequence of $(\lambda^n_\omega)$, we can find a further
subsequence $(\lambda^{n_k}_\omega)$ (which depends on $\omega$), a probability
space $(\esprob',\tribu',\prob')$, and random $\espE$-valued variables 
$X_1,\dots,X_k,\dots$ and $X$ defined on $\esprob'$ such that the law
of $X_k$ is $\lambda^{n_k}_\omega$ for each $k$, the law of $X$ is
$\mu_\omega$, and $(X_k)$ converges $\prob'$-a.e.~to $X$. For such an
$\omega$, we have, by the dominated convergence theorem, 
\begin{multline*}
\lim_k \int_\espE f(\omega,x)\,d\lambda^{n_k}_\omega(x)
=\lim_k \int_{\esprob'} f(\omega,X_k)\,d\prob'\\
= \int_{\esprob'} f(\omega,X)\,d\prob'
=\int_\espE f(\omega,x)\,d\mu_\omega(x).
\end{multline*}
Thus, for $\omega$  in the almost sure set of \eqref{eq:komlos}, every
subsequence of $(\lambda^{n_k}_\omega)$ has a further subsequence for
which 
the convergence in \eqref{eq:komlos-f} holds. 
This proves \eqref{eq:komlos-f}. 
We deduce that, for any subsequence of $(\mu^n)$ we can extract a
further subsequence $(\nu^n)$ such that 
\begin{equation}\label{eq:komlos-f-global}
\lim_n\frac{1}{n}\sum_{k=1}^{n}{\nu}^n(f)
=\mu(f),
\end{equation}
which proves the lemma. 
\finpr

The following technical lemma will be useful for limits of integrals
of unbounded integrands with respect to Young measures.
\begin{lem}\label{lem:quadraticYoung}
Let $\espE$ be a Suslin submetrizable topological space, and let 
$(X_n)$ be a sequence of $\espE$-valued random variables defined on
$\esprob$. 
Assume that $(X_n)$ stably converges 
to a Young measure $\mu\in\youngs(\esprob,\tribu,\prob;\espE)$ 
(where each $X_n$ is identified with the Young measure
$\delta_{X_n(\omega)}\otimes\,d\prob(\omega)$). 
Let $\Phi :\, \esprob\times\espE\rightarrow\R$ be measurable such that
\begin{enumerate}[(i)]
\item $\Phi(\omega,.)$ is
sequentially continuous for all $\omega\in\esprob$, 
\item The sequence $(\Phi(.,X_n))$ is uniformly integrable.
\end{enumerate}
Then $\Phi$ is $\mu$-integrable, and 
\begin{equation*}
  \label{eq:pyoung-converge}
  \lim_{n}\expect\Phi(.,X_n)
=\int_{\esprob\times\espE} \Phi\,d\mu.
\end{equation*}
\end{lem}
\proof
We only need to prove Lemma \ref{lem:quadraticYoung} in the case when
$\Phi\geq 0$, the general result comes from $\Phi=\Phi_+-\Phi_-$.

For each $N\geq 0$, 
we have 
\begin{equation}\label{eq:UI}
\lim_{N\rightarrow+\infty}\sup_{n}
   \expect\CCO{\Phi(.,X_n)\un{\accol{\Phi(.,X_n)\geq N}}}=0. 
\end{equation}
Set 
$$\Phi^N=\begin{cases}
\Phi & \text{ if }\Phi\leq N\\
N   & \text{ if }\Phi\geq N.
\end{cases}$$
From the definition of stable convergence and Lemma \ref{rem:seq-continu}, 
we have, for each $N$, 
\begin{equation}\label{eq:stable-tronq}
\lim_{n\rightarrow+\infty}\expect\Phi^N(.,X_n)=\mu(\Phi^N).
\end{equation}
Furthermore, by \eqref{eq:UI}, the convergence in
\eqref{eq:stable-tronq} is uniform with respect to $N$. We thus have,
with the help of Beppo Levi's lemma:
\begin{align*}
\mu(\Phi)
&=\sup_N\mu(\Phi^N)=\lim_N\mu(\Phi^N)\\
&=\lim_N\lim_{n}\expect\Phi^N(.,X_n)\\
&=\lim_n\lim_{N}\expect\Phi^N(.,X_n)\\
&=\lim_n\expect\Phi(.,X_n).
\end{align*}%
\finpr

\paragraph{Construction of the extended probability space: the
  processes $Y$, $\VV$ and $Z$}

Recall that $\VV^{(n)}=\int_0^. Z_s^{(n)}\,d\wien_s$.
By Proposition \ref{lem:Ytight}, 
the sequence
$(Y^{(n)},\VV^{(n)})$, 
seen as a sequence of random variables with values in 
$\trajYDs\times\trajYDs$, 
is tight. 
Let us denote $\trajZH={\ellp{2}_{\espZ}([0,T])}$, and let
$\trajZHw$ be the space $\trajZH$ endowed with its weak topology (note
that this topology
has the same Borel sets as the strong topology). Each
$Z^{(n)}$ can be considered as a random variable with values in
$\trajZHw$. Furthermore, by Corollary \ref{cor:Zbounded}, 
the sequence $(Z^{(n)})$ is tight in $\trajZHw$: 
Indeed, the closed
balls are compact in  $\trajZHw$, and we have 
\begin{align*}
\sup _n\prob\accol{\norm{Z^{(n)}}_{\trajZH}\geq R}
&\leq \sup_n \frac{1}{R^2\,}\expect\int_0^T \norm{Z^{(n)}_s}_{\trajZ}^2ds\\
&\rightarrow 0 \text{ when $R\rightarrow\infty$.}
\end{align*}
Thus $(Y^{(n)},\VV^{(n)},Z^{(n)})$ is a tight sequence of 
$\trajYDs\times\trajYDs\times\trajZHw$-valued variables.

We now consider the space
$\youngs(\esprob,\tribu,\prob;\trajYDs\times\trajYDs\times\trajZHw)$, 
which we denote for simplicity by $\youngs$.  
By 
Prohorov's sequential compactness criterion for Young measures 
\cite[Theorem 4.3.5]{cc-prf-valadier04book},
we can extract a subsequence of $(Y^{(n)},\VV^{(n)},Z^{(n)})$ 
(for simplicity, we denote this extracted sequence 
by $(Y^{(n)},\VV^{(n)},Z^{(n)})$)
which converges 
{stably} to some $\mu\in\youngs$, 
that is, for every measurable bounded mapping 
$\Phi : \, \esprob\times\trajYDs\times\trajYDs\times\trajZHw\mapsto \R$ 
such that $\Phi(\omega,.,.,.)$ is continuous for all $\omega$, 
we have
\begin{multline}
\lim_{n\rightarrow\infty}
\int_\esprob \Phi\CCO{\omega,Y^{(n)}(\omega),\VV^{(n)}(\omega),Z^{(n)}(\omega)}
                                                    \,d\prob(\omega)\\
=\int_\esprob\int_{\trajYD\times\trajYD\times\trajZH}\Phi(\omega,y,\vv,z)
                                     \,d\mu_\omega(y,\vv,z)\,d\prob(\omega). 
\label{eq:young-convergence}
\end{multline}
In particular, $(Y^{(n)},\VV^{(n)},Z^{(n)})$ converges in law to the image of
$\mu$ by the canonical projection of
 $\esprob\times\trajYDs\times\trajYDs\times\trajZHw$
on $\trajYDs\times\trajYDs\times\trajZHw$.

Let us denote by $\Ctribu$  
the Borel
 $\sigma$-algebra of $\trajYD$ (recall that S has the same Borel
 subsets as Skorokhod's $J_1$ topology), 
 and, for each $t\in [0,T]$, let $\Ctribu_t$ be the sub-$\sigma$-algebra
 of $\Ctribu$ generated by the projection onto $\trajYD_\espY([0,t])$. 
Similarly, let $\Htribu$  denote 
the Borel
 $\sigma$-algebra of $\trajZHw$, and, for each $t\in [0,T]$, 
 let $\Htribu_t$ be the sub-$\sigma$-algebra
 of $\Htribu$ generated by the projection onto ${\ellp{2}_{\espZ}([0,t])}$.  
We define a stochastic basis $(\esprobb,\tribuu,(\tribuu_t)_t,\mu)$ by 
\begin{equation*}
\esprobb=\esprob\times\trajYD\times\trajYD\times\trajZH,
\quad 
\tribuu=\tribu\otimes\Ctribu\otimes\Ctribu\otimes\Htribu,
\quad 
\tribuu_t=\tribu_t\otimes\Ctribu_t\otimes\Ctribu_t\otimes\Htribu_t,
\end{equation*}
and we define a process
$(Y,\ZK,Z)$ on $\esprobb$ by 
$$Y(\omega,\yy,\zk,z)=\yy, \quad \ZK(\omega,\yy,\zk,z)=\zk, 
\quad Z(\omega,\yy,\zk,z)=z.$$
Clearly, 
$(Y,\ZK,Z)$ is 
$(\tribuu_t)$--adapted. 
Furthermore, the law of $(Y,\ZK,Z)$  
is the marginal measure of $\mu$ on $\trajYD\times\trajYD\times\trajZH$, 
in particular 
$(Y^{(n)},\VV^{(n)},Z^{(n)})$ converges in law to $(Y,\ZK,Z)$ 
on $\trajYDs\times\trajYDs\times\trajZHw$.  
By \cite[Theorem 3.11]{jakubowski97skor}, 
{\em we can (and will) furthermore choose 
the extracted sequence such that, there exists a countable set 
$\countabl\subset[0,T[$ such that, for every $t\in[0,T]\setminus\countabl$, 
the sequence $(Y^{(n)}_t,\VV^{(n)}_t)$ converges in law to $(Y_t,\ZK_t).$}

Now, the random variables $(Y^{(n)},\VV^{(n)},Z^{(n)})$ 
can be seen as random elements defined
on $\esprobb$, using the notations, for $n\geq 1$:
\begin{align*}
Y^{(n)}(\omega,\yy,\zk,z)&:=Y^{(n)}(\omega),\\
\VV^{(n)}(\omega,\yy,\zk,z)&:=\VV^{(n)}(\omega),\\
Z^{(n)}(\omega,\yy,\zk,z)&:=Z^{(n)}(\omega).
\end{align*}
Furthermore, $(Y^{(n)},\VV^{(n)},Z^{(n)})$ is $(\tribuu_t)$--adapted for each $n$. 
Likewise, we set $\wien(\omega,\yy,\zk,z)=\wien(\omega)$.

\begin{lem}\label{lem:wien}
The process $\wien$ is an $(\tribuu_t)$--standard Brownian motion under the
probability $\mu$. 
\end{lem}
\proof
By Balder's result on K-convergence  
\cite{balder89Prohorov,balder90new}, which is
valid for Hausdorff spaces with metrizable compact subsets 
\cite[Lemma 4.5.4]{cc-prf-valadier04book},  
each subsequence of 
 $({Y^{(n)},\VV^{(n)},Z^{(n)}})$
contains a further subsequence 
$(Y^{(n_k)},\VV^{(n_k)},Z^{(n_k)})$
which K--converges to $\mu$, that is, for each subsequence 
$(Y^{(n'_k)},\VV^{(n'_k)},Z^{(n'_k)})$ of
$(Y^{(n_k)},\VV^{(n_k)},Z^{(n_k)})$, 
we have 
$$\lim_n\frac{1}{n}\sum_{k=1}^{n}\dirac{(Y^{(n'_k)}(\omega),\VV^{(n'_k)}(\omega),
Z^{(n'_k)}(\omega))}
=\mu_\omega\ \text{a.e.}$$
where $\dirac{(y,\vv,z)}$ denotes the Dirac measure on $(y,\vv,z)$ and the
limit is taken in the narrow convergence.  
This entails that, 
for every $B\in\Ctribu_t\otimes\Ctribu_t\otimes\Htribu_t$, 
the mapping $\omega\mapsto\mu_\omega(B)$
is $\tribu_t$--measurable. 
The result follows from Lemma \ref{lem:extended_martingale}. 
\finpr

\paragraph{Properties of the processes $Y$ and $V$}
\begin{lem}\label{lem:classemonotone}
Let $\vargenh$ 
and $\vargenk$ 
be $\espY$-valued random variables 
defined on $\esprobb$. 
Let $t\in[0,T]$.
In order that $\vargenh$ and $\vargenk$ 
have the same conditional expectation with respect to $\tribuu_t$, 
it is sufficient that 
\begin{multline}
  \label{eq:classmonoton}
  \int_\esprob\int_{\trajYD\times\trajYD\times\trajZH}
     \Phi(\omega,y,\vv,z)\vargenh(\omega,y,\vv,z)
             \,d\mu_\omega(y,\vv,z)\,d\prob(\omega)\\
  =\int_\esprob\int_{\trajYD\times\trajYD\times\trajZH}
     \Phi(\omega,y,\vv,z)\vargenk(\omega,y,\vv,z)
             \,d\mu_\omega(y,\vv,z)\,d\prob(\omega)
\end{multline}
for every bounded $\tribuu_t$-measurable function
$\Phi :\esprobb\rightarrow\R$ such that 
$\Phi(\omega, .,.,.)$ is continuous for all $\omega\in\esprob$. 
\end{lem}
\proof
Let $\classA$ be the set of functions 
$\Phi :\esprobb\rightarrow\R$ which are $\tribuu_t$-measurable
and such that 
$\Phi(\omega, .,.,.)$ is continuous for all $\omega\in\esprob$. 
The set $\classA$
is stable by multiplication of two functions and generates $\tribuu_t$. 
Assume that \eqref{eq:classmonoton} holds for every 
$\Phi\in\classA$, and 
let $\classB$ be the vector space of bounded $\tribuu_t$-measurable 
functions $\varphi$ defined on $\esprobb$ such that
\begin{multline*}
  \int_\esprob\int_{\trajYD\times\trajYD\times\trajZH}
     \varphi(\omega,y,\vv,z)\vargenh(\omega,y,\vv,z)
             \,d\mu_\omega(y,\vv,z)\,d\prob(\omega)\\
  =\int_\esprob\int_{\trajYD\times\trajYD\times\trajZH}
     \varphi(\omega,y,\vv,z)\vargenk(\omega,y,\vv,z)
             \,d\mu_\omega(y,\vv,z)\,d\prob(\omega).
\end{multline*}
The space $\classB$ contains $\classA$. Furthermore, 
$\classB$ contains the constant functions and 
is stable under monotone limits of uniformly
bounded sequences. 
By the  monotone class theorem
(see 
\cite[Appendix A0]{sharpe88book} and 
 \cite[Th´eor\`eme 21, page 20]{dellacheriemeyer75book}), 
$\classB$ contains all bounded $\tribuu_t$-measurable functions. 
\finpr

\begin{lem}\label{lem:Vmartingale}
The process $\VV$ is a martingale with respect to
$(\esprobb,\tribuu,(\tribuu_t)_t,\mu)$. 
\end{lem}
\proof
Let $t\in[0,T]$, and let $s\in[0,T-t]$. 
By Lemma \ref{lem:classemonotone}, 
in order to prove that $\espcond{\VV_{t+s}}{\tribuu_t}=\VV_t$, 
we only need to show that, 
for each bounded $\tribuu_t$-measurable 
$\Phi :\esprobb\rightarrow\R$ such that 
$\Phi(\omega, .,.,.)$ is continuous for all $\omega\in\esprob$, 
we have
\begin{equation}\label{eq:VVmg}
\expect\CCO{ \Phi \times\VV_{t+s}}
=\expect\CCO{\Phi \times\VV_{t}}.
\end{equation}
Let us denote, 
for any  $r\in[0,T]$, any $\vv\in\trajYD$ 
and any $\delta>0$
\begin{equation}\label{eq:proj}
\proj{r}(\vv)=\vv(r)\ \text{ and }\ 
\proj{r,\delta}(\vv)=\frac{1}{\delta}\int_r^{r+\delta}\vv(s)\,ds.
\end{equation} 
The mapping $\proj{r} :\,\trajYD\rightarrow \espY$ is not continuous
for the topology S,
but $\proj{r,\delta}$ is S-continuous, and we have 
$$\lim_{\delta\rightarrow0}\proj{r,\delta}(\vv)=\proj{r}(\vv).$$  
Let $\delta>0$. 
Let 
$$\phi(\omega,y,\vv,z)
=\Phi(\omega,y,\vv,z)\CCO{\proj{t+s,\delta}(\vv)-\proj{t,\delta}(\vv)}.
$$ 
By Corollary \ref{cor:Zbounded}, the sequence 
$(\phi(\omega,Y^{(n)},\VV^{(n)},Z^{(n)}))$ 
is bounded in $\ellp{2}_{\espY}(\esprob)$, thus it is uniformly
integrable. We can thus apply  Lemma \ref{lem:quadraticYoung} to the
integrand $\phi$. 
Using the definition of $\VV$  
and the fact that each $\VV^{(n)}$ is a martingale, we get
\begin{multline*}
\expectmu\CCO{\Phi \times
   \CCO{\frac{1}{\delta}\int_{t+s}^{t+s+\delta}\VV_{u}\,du
        -\frac{1}{\delta}\int_{t}^{t+\delta}\VV_u\,du}}\\
\begin{aligned}
&=\int_\esprob\int_{\trajYD\times\trajYD\times\trajZH}
      \Phi(\omega,y,\vv,z)\CCO{\proj{t+s,\delta}-\proj{t,\delta}}(v)
              \,d\mu_\omega\CCO{y,\vv,z}\,d\prob(\omega)\\
&=\lim_{n\rightarrow\infty}
   \int_\esprob
   \Phi\CCO{\omega,Y^{(n)}(\omega),\VV^{(n)}(\omega),Z^{(n)}(\omega)}
        \frac{1}{\delta}\int_{t}^{t+\delta}
                  \CCO{\VV^{(n)}_{u+s}(\omega)-\VV^{(n)}_{u}(\omega)}\,du
                                              \,d\prob(\omega)\\
&=\lim_{n\rightarrow\infty}
   \int_\esprob
      \Phi\CCO{\omega,Y^{(n)}(\omega),\VV^{(n)}(\omega),Z^{(n)}(\omega)}
        \espcondF{\frac{1}{\delta}\int_{t}^{t+\delta}
                \CCO{\VV^{(n)}_{u+s}-\VV^{(n)}_{u}}\,du}{t}
                                              \,d\prob\\
&=0.
\end{aligned}
\end{multline*}
We deduce that 
\begin{multline*}
\expectmu\CCO{\Phi \times\CCO{\VV_{t+s}-\VV_t}}\\
=\lim_{\delta\rightarrow 0}\expectmu\CCO{\Phi \times
   \CCO{\frac{1}{\delta}\int_{t+s}^{t+s+\delta}\VV_{u}\,du
        -\frac{1}{\delta}\int_{t}^{t+\delta}\VV_u\,du}}
=0
\end{multline*}
by boundedness in $\ellp{2}_\espY(\esprobb,\tribuu,\mu)$ of
$(\VV_r)_{0\leq r\leqq T}$.
\finpr

\begin{lem}\label{lem:LMortogonaux}
Let $\MZ_t=\int_0^t Z_s\,d\wien_s$.
The martingale
$L:=\VV-\MZ$ is orthogonal to $W$. 
\end{lem}
\proof
Let us denote the coordinates processes as in the following examples: 
$\VV=(\VV^\coord{i})_{1\leq i\leq d}$,
$Z^{(n)}_t=(Z^{(n),\coord{i,k}}_t)_{1\leq i\leq d,1\leq k\leq m}$,  
$Z_t=(Z^{\coord{i,k}}_t)_{1\leq i\leq d,1\leq k\leq m}$,
$\wien_t=(\wien^{\coord{k}})_{1\leq k\leq m}$. 

Let $i\in\{1,\dots,m\}$ and $j\in\{1,\dots,d\}$. 
Let us denote by $\cvarqd{P}{Q}$ the quadratic cross variation of two
semimartingales $P$ and $Q$. For each $n$, let
\begin{align*}
\mgx^{(n),\coord{i,j}}_t
&=\wien^{\coord{i}}_{t}{\VV^{(n),\coord{j}}_{t}}
  -\cvarqd{\wien^{\coord{i}}}{\VV^{(n),\coord{j}}}_t\\
&=\wien^{\coord{i}}_t
\sum_{k=1}^m \int_0^t {Z^{(n),\coord{j,k}}_r}\,d\wien^{\coord{k}}_r
-\int_0^t{Z^{(n),\coord{j,i}}_r}\,dr\\ 
\mgxx^{\coord{i,j}}_t
&=\wien^{\coord{i}}_{t}{\MZ^{\coord{j}}_{t}}
-\cvarqd{\wien^{\coord{i}}}{\MZ^{\coord{j}}}_t\\
&=\wien^{\coord{i}}_t
\sum_{k=1}^m\int_0^t {Z^{\coord{j,k}}_r}\,d\wien^{\coord{k}}_r
-\int_0^t {Z^{\coord{j,i}}_r}\,dr.
\end{align*}
As $\wien$ is continuous, the processes $\mgx^{(n),\coord{i,j}}$ and 
$\mgxx^{\coord{i,j}}$ are continuous martingales. 
Let $\Phi :\esprobb\rightarrow\R$ be a bounded $\tribuu_t$-measurable 
function such that 
$\Phi(\omega, .,.,.)$ is continuous for all $\omega\in\esprob$. 
Observe that,
from the 
stable convergence 
of $\VV^{(n)}$ to $\VV$, 
we have, for any $\tau\in[0,T]$ and any $\delta>0$, 
\begin{multline}
\lim_n\expect\CCO{
  { \proj{\tau,\delta}(\wien^{\coord{i}}\VV^{(n),\coord{j}}) }
                                        \Phi\CCO{.,Y^{(n)},\VV^{(n)},Z^{(n)}}} \\
=\expect\CCO{
    { \proj{\tau,\delta}(\wien^{\coord{i}}\VV^{\coord{j}})}\Phi
             }
\label{eq:MVn}
\end{multline}
using {Lemma \ref{lem:quadraticYoung}} with the integrand 
$\phi(\omega,y,\vv,z)
   =\proj{\tau,\delta}(\wien^{\coord{i}}(\omega)\vv^{(n),\coord{j}})
                       \Phi(\omega,y,\vv,z)$,
         where $\proj{\tau,\delta}$ is defined as in \eqref{eq:proj}.
Similarly, 
from the 
stable convergence 
of $(\VV^{(n)},Z^{(n)})$ to $(\VV,Z)$, and 
applying {Lemma \ref{lem:quadraticYoung}} with the integrand 
$$\phi(\omega,y,\vv,z)=
\CCO{
\int_0^{\tau} {z^{\coord{j,i}}_r}\,dr}
\Phi(\omega,y,\vv,z),             
$$
 we get
\begin{multline}
\label{eq:HZn}
\lim_n\expect\CCO{
    {\left\lgroup
\int_0^{\tau}
        {Z^{(n),\coord{j,i}}_r}\,dr
    \right\rgroup}
\,\Phi\CCO{.,Y^{(n)},\VV^{(n)},Z^{(n)}}
}\\
=\expect\CCO{
    \left\lgroup{
\int_0^{\tau}
        {Z^{\coord{j,i}}_r}\,dr
    }\right\rgroup
\,\Phi
}.
\end{multline}
Let $t\in[0,T]$, and let $s\in[0,T-t]$. 
Using \eqref{eq:MVn}, \eqref{eq:HZn}, and 
the fact that $\mgx^{(n),\coord{i,j}}$ and $\mgxx^{\coord{i,j}}$ are
martingales, 
we get, for any $\delta>0$,
\begin{align*}
&\frac{1}{\delta}\int_{t}^{t+\delta}\expect\CCO{
    \wien^{\coord{i}}_{u+s}\CCO{\VV^{\coord{j}}_{u+s}-\MZ^{\coord{j}}_{u+s}}
                                                \Phi }du\\
=&\frac{1}{\delta}\int_{t}^{t+\delta}\lim_n\Biggl\lgroup\expect\CCO{
   \wien^{\coord{i}}_{u+s} { \VV^{(n),\coord{j}}_{u+s} }   
              \Phi\CCO{.,Y^{(n)},\VV^{(n)},Z^{(n)}} }\\ 
&\phantom{\frac{1}{\delta}\int_{t}^{t+\delta}}
  -\expect\CCO{
    \wien^{\coord{i}}_{u+s} { \MZ^{\coord{j}}_{u+s} }
    \Phi\CCO{.,Y^{(n)},\VV^{(n)},Z^{(n)}} }\Biggr\rgroup du\\
=&\frac{1}{\delta}\int_{t}^{t+\delta}
 \lim_n\Biggl\lgroup\expect\CCO{
\CCO{\mgx^{(n),\coord{i,j}}_{u+s} 
+\int_0^{u+s} {Z^{(n),\coord{j,i}}_r}\,dr
    }
\Phi\CCO{.,Y^{(n)},\VV^{(n)},Z^{(n)}}
}\\
&\phantom{\frac{1}{\delta}\int_{t}^{t+\delta}}
-\expect\CCO{
\CCO{\mgxx^{\coord{i,j}}_{u+s} 
+\int_0^{u+s}{Z^{\coord{j,i}}_r}\,dr
    }
\Phi\CCO{.,Y^{(n)},\VV^{(n)},Z^{(n)}}
}\Biggr\rgroup du\\
=&\frac{1}{\delta}\int_{t}^{t+\delta}
\lim_n\expect\CCO{
{\CCO{\mgx^{(n),\coord{i,j}}_{u+s}-\mgxx^{\coord{i,j}}_{u+s}}
    }
\Phi\CCO{.,Y^{(n)},\VV^{(n)},Z^{(n)}}
}du
\\
=&\frac{1}{\delta}\int_{t}^{t+\delta}
\lim_n\expect\CCO{
{\CCO{\mgx^{(n),\coord{i,j}}_{u}-\mgxx^{\coord{i,j}}_{u}}
    }
\Phi\CCO{.,Y^{(n)},\VV^{(n)},Z^{(n)}}
}du
\\
=&\frac{1}{\delta}\int_{t}^{t+\delta}
\expect\CCO{
    \wien^{\coord{i}}_{u}\CCO{\VV^{\coord{j}}_{u}-\MZ^{\coord{j}}_{u}}
                                                \Phi }du.
\end{align*}  
Passing to the limit when $\delta\rightarrow0$ yields
$$
\expect\CCO{
    \wien^{\coord{i}}_{t+s}\CCO{\VV^{\coord{j}}_{t+s}-\MZ^{\coord{j}}_{t+s}}
                                                \Phi }
=\expect\CCO{
    \wien^{\coord{i}}_{t}\CCO{\VV^{\coord{j}}_{t}-\MZ^{\coord{j}}_{t}}
                                                \Phi}.
$$
By Lemma \ref{lem:classemonotone}, this shows that 
$\wien^{\coord{i}} \CCO{ \VV^{\coord{j}}-\MZ^{\coord{j}} }$ 
is a martingale.

\finpr

\subsection{Proof of the main result}\label{subsect:proof-main}
In this part, we use the special form of $f$ with respect to $Z$: 
By hypothesis (\Hbb), $f$ has the form
\begin{equation}\label{eq:affine}
f(s,x,y,z)=\alpha(s,x,y)z+\beta(s,x,y),
\end{equation} 
where $\alpha$ and $\beta$
are bounded and continuous in $(x,y)$, and $\alpha$ takes its
values in the space $\linear(\espZ,\esp)$  
of linear mappings from $\espZ$ to $\esp$. 

We first prove a technical lemma.

\begin{lem}\label{lem:joint_continuity}
Let $\esplin$ be the space of linear mappings from $\espY$ to $\R^{l}$ 
for some $l\geq 1$. 
Let $b :\,[0,T]\rightarrow\esplin$ be a continuous function. 
For each $t\in[0,T]$,
the mapping 
$$\Phi :\,
\left\{\begin{array}{lcl}
\trajYDs\times\trajZHw&\rightarrow&\R^{l}\\
(y,z)&\mapsto&\int_0^t b(s).f(s,x(s),y(s),z(s))\,ds
\end{array}\right.
$$
is sequentially continuous. 
Furthermore, if $y_n\rightarrow y$ in $\trajYDs$ and $z_n\rightarrow z$
in $\trajZHw$, then, for every $t\in[0,T]$, we have  
\begin{multline*}
\lim_n\biggl(
 \int_0^t b(s).f(s,x(s),y_n(s),z_n(s+1/n))\,ds\\
- \int_0^t b(s).f(s,x(s),y_n(s),z_n(s))\,ds\biggr)
=0.
\end{multline*}
\end{lem}
\proof
We only need to prove the lemma for 
$f(s,x,y,z)=\alpha(s,x,y)z$. 
As $x$ does not play any role in our
reasoning, we write for simplicity 
$f(s,x,y,z)=\alpha(s,y)z$. 

First, for every $\zz\in\trajZ$, we have 
\begin{equation}
  \label{eq:znl2}
  \lim_{n\rightarrow\infty}\norm{z-z(.+1/n)}_{\trajZ}=0.
\end{equation}
Indeed, for every $\epsilon>0$, 
there exists a continuous function $u :\,[0,T]\rightarrow\espZ$ such
that $\norm{z-u}_{\trajZ}<\epsilon$. Then we have, for every $n\geq
1$, 
$$\norm{z(.+1/n)-u(.+1/n)}_{\trajZ}<\epsilon.$$
But the family $u(.+1/n)$ is uniformly integrable because it is
bounded in $\trajZ$, thus, 
by Vitali's theorem and the continuity of $u$, 
$$\lim_{n\rightarrow\infty}\norm{u-u(.+1/n)}_{\trajZ}=0.$$
We conclude by the triangular inequality that 
$$\limsup_{n\rightarrow\infty}\norm{z-z(.+1/n)}_{\trajZ}\leq 2\epsilon,$$
which proves \eqref{eq:znl2}.

Now, 
let $y_n\rightarrow y$ in $\trajYDs$ and $z_n\rightarrow z$
in $\trajZHw$. 
We have in particular
$$y_n(s)\rightarrow y(s)\text{ for a.e.~}s\in[0,T]
\text{ and }
\sup_n\norm{z_n}_{\trajZH}<+\infty,$$
thus
\begin{multline*}
\norm{\Phi(y_n,z_n)(t)-\Phi(y,z)(t)}\\
\begin{aligned}
=&\biggl\Vert\int_0^tb(s).\CCO{\alpha(s,y_n(s))-\alpha(s,y(s))}z_n(s)\,ds\\
&+\int_0^t b(s).\alpha(s,y(s))\CCO{z_n(s)-z(s)}ds\biggr\Vert\\
\leq & \sup_n \norm{z_n}_{\trajZH}
               \,\CCO{\int_0^t
                 \norm{b(s)}\norm{\alpha(s,y_n(s))-\alpha(s,y(s))}^2ds}^{1/2}\\
&+\norm{\int_0^t b(s).\alpha(s,y(s))\CCO{z_n(s)-z(s)}ds}\\
\rightarrow& \,0 \text{ when }n\rightarrow\infty,
\end{aligned}
\end{multline*}
which proves the first part of Lemma
\ref{lem:joint_continuity}. Furthermore, we have 
\begin{multline*}
\biggl\Vert
\int_0^t b(s).\alpha(s,y_n(s))\bigl({z_n(s+1/n)-z_n(s)}\bigr)\,ds\biggr\Vert
\\
\begin{aligned}
=&\biggl\Vert
  \int_0^t b(s).\bigl({\alpha(s,y_n(s))-\alpha(s,y(s))}\bigr)
   \bigl({z_n(s+1/n)-z_n(s)}\bigr)\,ds\\
&+\int_0^t  b(s).\alpha(s,y(s)) \bigl({z_n(s+1/n)-z_n(s)}\bigr)\,ds
\biggr\Vert\\
\leq&2\sup_n \norm{z_n}_{\trajZH}
      \,\CCO{\int_0^t
       \norm{b(s)}\norm{ \alpha(s,y_n(s))-\alpha(s,y(s)) }^2\,ds 
             }^{1/2}\\
&+\biggl\Vert
\int_0^t  b(s).\alpha(s,y(s)) \bigl({z_n(s+1/n)-z_n(s)}\bigr)\,ds
\biggr\Vert.
\end{aligned}
\end{multline*}
The term $\int_0^t\norm{b(s)}\norm{ \alpha(s,y_n(s))-\alpha(s,y(s)) }^2\,ds$
converges to 0 by the dominated convergence theorem. 
On the other hand, since $b$ and $\alpha$ are bounded 
and  $(z_n)$ is uniformly bounded in $\trajZ$,
we have (with the convention that $b(s)=\alpha(s)=0$ for $s<0$):
\begin{multline*}
\lim_n\,\biggl\Vert
\int_0^t b(s).\alpha(s,y(s))\bigl({z_n(s+1/n)-z_n(s)}\bigr)\,ds
\biggr\Vert
\\
\begin{aligned}
=&\lim_n\,\biggl\Vert
       \int_0^t  \bigl(b(s-1/n).\alpha(s-1/n,y(s-1/n))
        -b(s).\alpha(s,y(s))z_n(s)\bigr)\,ds
\biggr\Vert.
\end{aligned}
\end{multline*}
This term
vanishes
by \eqref{eq:znl2} with $z(s)=b(s).\alpha(s,y(s))$, using again 
the uniform boundedness of $(z_n)$ in $\trajZ$.
Thus 
$$\lim_n\,
\biggl\Vert
\int_0^t  b(s).\alpha(s,y_n(s)) \bigl({z_n(s+1/n)-z_n(s)}\bigr)\,ds
\biggr\Vert=0.$$
\finpr

In order to check that $(Y,Z)$ is a solution to \eqref{eq:BSDE-gene},
we prove in the next lemma that 
we can replace $\widetilde{Z}^{(n)}$ by $Z^{(n)}$ in the limit of 
$\int_t^Tf(s,X_s,Y_s^{(n)},\widetilde{Z}_s^{(n)})\,ds$.

\begin{lem}\label{lem:equivalence}
For each $t\in[0,T]$, the sequence 
$${\int_t^Tf(s,X_s,Y_s^{(n)},\widetilde{Z}_s^{(n)})\,ds
-
\int_t^Tf(s,X_s,Y_s^{(n)},{Z}_s^{(n)})\,ds
}$$
converges to 0 in $\mu$-probability. 
\end{lem}
\proof
With the notations of \eqref{eq:affine}, 
we only need to check that 
$${\int_t^T\alpha(s,X_s,Y_s^{(n)})\bigl(\widetilde{Z}_s^{(n)}
                                        -{Z}_s^{(n)}\bigr)\,ds
}$$
converges to 0 in probability. 

Now, by Lemma \ref{lem:Zbounded} and Proposition \ref{prop:YZbounded}, 
the sequence $(\widetilde{Z}^{(n)})$ is bounded in 
${\ellp{2}_{\espZ}(\esprob\times[0,T])}$, thus it can be viewed as a
tight sequence of $\trajZHw$-valued random variables. 
Enlarging the space
$\esprobb$ to
$\esprob\times\trajYD\times\trajYD\times\trajZH\times\trajZH$, we can 
assume that $(Y^{(n)},\VV^{(n)},Z^{(n)},\widetilde{Z}^{(n)})$
converges 
to a Young measure, still  denoted by $\mu$, 
in 
$\youngs(\esprob,\tribu,\prob;
                   \trajYDs\times\trajYDs\times\trajZHw\times\trajZHw)$. 
We set 
$$\widetilde{Z}(\omega,\yy,\zk,z,\widetilde{z})=\widetilde{z}$$ 
and we
extend $Y$, $\VV$, $Z$, and the $\sigma$-algebra $\tribuu_t$ 
in the obvious way.  

Let $\vargenk$ be an $\tribuu_t$-adapted process with \cadlag\
trajectories 
in $\espY$, and assume
that  $\vargenk$ is continuous with respect to $y$, $\vv$, $z$, and
$\widetilde{z}$, 
and that  
the sequence 
$$\bigl(\int_0^T \vargenk_s
.\alpha(s,X_s,Y_s)(Z_s^{(n)}-\widetilde{Z}_s^{(n)})\,ds\bigr)$$ 
is uniformly integrable. 
We have
\begin{multline*}
\expect\int_0^T \vargenk_s
.\alpha(s,X_s,Y_s)(Z_s-\widetilde{Z}_s)\,ds\\
=\lim_n  
\expect\int_0^T \vargenk_s .\alpha(s,X_s,Y^{(n)}_s)
                     (Z_s^{(n)}-\widetilde{Z}^{(n)}_s)\,ds
\end{multline*}
by 
Lemma \ref{lem:joint_continuity} and 
Lemma \ref{lem:quadraticYoung}, with 
$$\Phi(\omega,y,\vv,z,\widetilde{z})
=\int_0^T
\vargenk_s(\omega,y,\vv,z,\widetilde{z}).
       \alpha(s,x_s,y_s)(z_s-\widetilde{z}_s)\,ds.$$
Thus, by Lemma \ref{lem:joint_continuity},
\begin{multline*}
\expect\int_0^T \vargenk_s
.\alpha(s,X_s,Y_s)(Z_s-\widetilde{Z}_s)\,ds\\
\begin{aligned}
=&\lim_n  
\expect\int_0^T \vargenk_s .\alpha(s,X_s,Y^{(n)}_s)
                     (Z_{s+1/n}^{(n)}-\widetilde{Z}^{(n)}_s)\,ds\\
=&\lim_n  
\expect\int_0^T \espcondFbar{\vargenk_s .\alpha(s,X_s,Y^{(n)}_s)
                     (Z_{s+1/n}^{(n)}-\widetilde{Z}^{(n)}_s)}{s}\,ds\\
=&\lim_n  
\expect\int_0^T \vargenk_s .\alpha(s,X_s,Y^{(n)}_s)
                     (\widetilde{Z}^{(n)}_s-\widetilde{Z}^{(n)}_s)\,ds\\
=&0.
\end{aligned}
\end{multline*}
In particular, one can take 
$$K_s=\frac{\alpha(s,X_s,Y_s)(Z_s-\widetilde{Z}_s)}%
{1+{\CCO{\alpha(s,X_s,Y_s)(Z_s-\widetilde{Z}_s)}}^2}.$$
Thus $\alpha(s,X_s,Y_s)(Z_s-\widetilde{Z}_s)=0$, $\mu$-a.e.,~for almost
every $s\in[0,T]$.

Let $\Psi :\,\espY\rightarrow\R$ be a bounded continuous function. 
Let 
$$\Phi :\,\left\{\begin{array}{lcl}
 \trajYDs\times\trajYDs\times\trajZHw\times\trajZHw&\rightarrow&\R \\
  (x,y,z,\widetilde{z})&\mapsto&
             \Psi\CCO{\int_{t}^T \alpha(s,x_s,y_s)({z}_s-\widetilde{z}_s)\,ds}.
\end{array}\right.$$ 
By Hypothesis (\Hbb) and 
Lemma \ref{lem:joint_continuity}, 
$\Phi$ is sequentially continuous, 
thus, from 
the $\tribu$-stable convergence of 
$(X,Y^{(n)},{Z}^{(n)},\widetilde{Z}^{(n)})$,
\begin{multline*}
\Psi(0)=
\expect\Psi\CCO{\int_{t}^T\alpha(s,X_s,Y_s)({Z}_s-\widetilde{Z}_s)\,ds}
=\lim_n \expect \Phi(X,Y^{(n)},{Z}^{(n)},\widetilde{Z}^{(n)})\\
=\lim_n \expect\Psi\CCO{\int_{t}^T
     \alpha(s,X_s,Y^{(n)}_s)\bigl(\widetilde{Z}_s^{(n)}
                                        -{Z}_s^{(n)}\bigr)\,ds}.
\end{multline*}
This shows that the sequence 
$$\CCO{\int_t^T\alpha(s,X_s,Y_s^{(n)})\bigl(\widetilde{Z}_s^{(n)}
                                        -{Z}_s^{(n)}\bigr)\,ds
}$$
converges to 0 in law, thus in probability. 
\finpr

\begin{lem}\label{lem:bicontinuous}
The sequence  
$(\int_{.}^Tf(s,X_s,Y_s^{(n)},\widetilde{Z}_s^{(n)})\,ds)$ 
converges in law to 
$\int_{.}^Tf(s,X_s,Y_s,{Z}_s)\,ds$. 
\end{lem}
\proof
By Lemma \ref{lem:ftight},  we know that the sequence 
  $(\int_{.}^Tf(s,X_s,Y_s^{(n)},\widetilde{Z}_s^{(n)})\,ds)$ 
is relatively compact in law,
 thus we only need to show that it has only one possible limit in law, 
and that this limit is 
the law  of $\int_{.}^Tf(s,X_s,Y_s,{Z}_s)\,ds$. 
By Lemma \ref{lem:equivalence}, it suffices to prove that 
$(\int_{t}^Tf(s,X_s,Y_s^{(n)},{Z}_s^{(n)})\,ds)$ converges
in law to 
$\int_{t}^Tf(s,X_s,Y_s,{Z}_s)\,ds$ for each $t\in[0,T]$. 

Let $\Psi :\,\espY\rightarrow\R$ be a bounded continuous function. 
Let 
$$\Phi :\,\left\{\begin{array}{lcl}
 \trajYDs\times\trajYDs\times\trajZHw&\rightarrow& \R \\
  (x,y,z)&\mapsto&\Psi\CCO{\int_{t}^T f(s,x_s,y_s,{z}_s)\,ds}.
\end{array}\right.$$ 
By Hypothesis (\Hbb) and Lemma
\ref{lem:joint_continuity}, 
$\Phi$ is sequentially continuous, 
thus, from
the $\tribu$-stable convergence of 
$(X,Y^{(n)},{Z}^{(n)})$,
\begin{multline*}
\expect\Psi\CCO{\int_{t}^Tf(s,X_s,Y_s,{Z}_s)\,ds}
=\mu(\Phi)
=\lim_n \expect \Phi(X,Y^{(n)},{Z}^{(n)})\\
=\lim_n \expect\Psi\CCO{\int_{t}^Tf(s,X_s,Y^{(n)}_s,{Z}^{(n)}_s)\,ds}.
\end{multline*}
\finpr

\proofof{Theorem \ref{theo:main}}
By Lemma \ref{lem:wien}, $\wien$ is a Brownian motion on
$(\esprobb,\tribuu,(\tribuu_t)_t,\mu)$.
Let $L_t=\VV_t-\VV_0-\MZ_t$, $0\leq t\leq T$. 
We have $L_0=0$ and  
$L$ is a \cadlag\ martingale by Lemma \ref{lem:Vmartingale}, 
furthermore   
 $L$ is orthogonal to $\wien$ by Lemma \ref{lem:LMortogonaux}. 
Thus 
there only remains to prove that 
$(Y,Z,L)$ satisfies \eqref{eq:BSDE-L}.

Thanks to Proposition \ref{lem:Ytight} and Lemma \ref{lem:ftight}, 
we know that the sequence 
\begin{equation}\label{eq:bigsequence}
\biggl( X,
        Y^{(n)}, 
        \int_{.}^Tf(s,X_s,Y_s^{(n)},\widetilde{Z}_s^{(n)})\,ds,
        \int_{.}^TZ_s^{(n)}\,d\wien_s\biggr)_{n\geq 1}
\end{equation}
is tight in
$\trajX\times\trajYDs\times\trajY\times\trajYDs$. Furthermore, 
$\biggl(\int_{.}^TZ_s^{(n)}\,d\wien_s\biggr)_{n\geq 1}$ converges in
law to $\VV_T-\VV_.$, and, by  
Lemma \ref{lem:bicontinuous}, 
$\biggl
(\int_{.}^Tf(s,X_s,Y_s^{(n)},\widetilde{Z}_s^{(n)})\,ds
\biggr)_{n\geq 1}$ converges in law to 
$\int_{.}^Tf(s,X_s,Y_s,Z_s)\,ds$. 
Extracting if necessary a further subsequence, 
we can thus assume that 
the sequence \eqref{eq:bigsequence}
jointly converges in law on
$\trajX\times\trajYDs\times\trajY\times\trajYDs$ to 
$$\biggl( X,
         Y, 
         \int_{.}^Tf(s,X_s,Y_s,Z_s)\,ds,
         \VV_T-\VV_.\biggr).$$
Then 
the process
$$U_.^{(n)}=Y^{(n)}_.-\xi-\int_{.}^T f(s,X_s,Y_s^{(n)},{Z}_s^{(n)})\,ds
        +\int_.^T Z_s^{(n)}\,d\wien_s 
        $$
converges in law in $\trajYDs$ to 
\begin{align*}
U_.:=&Y_.-\xi-\int_{.}^T f(s,X_s,Y_s,{Z}_s)\,ds
        +\VV_T-\VV_.\\
=&Y_.-\xi-\int_{.}^T f(s,X_s,Y_s,{Z}_s)\,ds
        +\int_.^TZ_s\,d\wien_s+L_T-L_.
\end{align*} 
But, by Lemma \ref{lem:Utn},  $(\sup_{0\leq t\leq T}U_t^{(n)})$ converges to 0 in
probability, thus $U=0$ a.e., which proves Theorem \ref{theo:main}. 
\finpr

\paragraph{Acknowledgements} We thank Adam Jakubowski
for introducing us to the topology S. %

We are also greatly indebted to the referees, who detected several
errors and gaps in previous versions, and helped improve 
the correct parts. 
We thank them for their thorough reading, 
their contributions, 
and their patience.


\end{document}